\documentclass[leqno, 10pt]{amsart}
\usepackage{amsmath, amsthm, amssymb, latexsym}
\newtheorem{theorem}{Theorem}[section]
\newtheorem{proposition}[theorem]{Proposition}
\newtheorem{lemma}[theorem]{Lemma}

\newtheorem{definition}[theorem]{Definition}

\theoremstyle{remark}

\newtheorem{remark}[theorem]{Remark}

  \newtheorem*{acknowledgements}{Acknowledgements}

\numberwithin{equation}{section}

\newcommand{\op}[1]{\operatorname{#1}}

\newcommand{\brak}[1]{\ensuremath{\langle\! #1\!\rangle}}

\newcommand{\Tr}{\ensuremath{\op{Tr}}}
\newcommand{\tr}{\op{tr}}

\newcommand{\Res}{\ensuremath{\op{Res}}}

\newcommand{\Str}{\op{Str}}

\newcommand{\C}{\ensuremath{\mathbb{C}}} 
\newcommand{\bH}{\ensuremath{\mathbb{H}}} 
\newcommand{\N}{\ensuremath{\mathbb{N}}} 
\newcommand{\R}{\ensuremath{\mathbb{R}}} 
\newcommand{\Z}{\ensuremath{\mathbb{Z}}}

\newcommand{\Rd}{\ensuremath{\R^{d+1}}}

\newcommand{\Rdo}{\R^{d+1}\!\setminus\! 0}

\newcommand{\URd}{U\times\R^{d+1}}

\newcommand{\URdo}{U\times(\R^{d+1}\!\setminus\! 0)}

\newcommand{\Ca}[1]{\ensuremath{\mathcal{#1}}}

\newcommand{\cE}{\Ca{E}}
\newcommand{\cF}{\Ca{F}}

\newcommand{\cI}{\ensuremath{\mathcal{I}}}

\newcommand{\cL}{\ensuremath{\mathcal{L}}}

\newcommand{\cN}{\ensuremath{\mathcal{N}}}

\newcommand{\fg}{\ensuremath{\mathfrak{g}}}
\newcommand{\fh}{\ensuremath{\mathfrak{h}}}

\newcommand{\psivdo}{$\Psi_{H}$DO}
\newcommand{\psivdos}{$\Psi_{H}$DO's}

\newcommand{\pvdo}{\ensuremath{\Psi_{H}}}

\newcommand{\pvdoz}{\ensuremath{\Psi_{H}^{\Z}}} 
 
\newcommand{\psido}{$\Psi$DO} 
\newcommand{\psidos}{$\Psi$DO's} 
\newcommand{\psinf}{\ensuremath{\Psi^{-\infty}}}

\newcommand{\pr}{\op{pr}}

\newcommand{\im}{\op{im}}

\newcommand{\End}{\ensuremath{\op{End}}}

\renewcommand{\Box}{\square}

\newcommand{\Sp}{\op{Sp}}

\newcommand{\dbarb}{\overline{\partial}_{b}}
 \newcommand{\dbarbE}{\overline{\partial}_{b,\cE}}
 \newcommand{\dbarbEq}{\overline{\partial}_{b,\cE;q}}
\newcommand{\dbarbqt}{\overline{\partial}_{b,t;q}}
\newcommand{\dbarbq}{\overline{\partial}_{b;q}}
\newcommand{\varthetabq}{\overline{\partial}^{*}_{b;q}}
\newcommand{\varthetabqt}{\overline{\partial}^{*}_{b,t;q}}
\newcommand{\Boxbqt}{\Box_{b,t;q}}

\begin{document}
\title{Noncommutative residue invariants\\
for CR and contact manifolds} 

\author{Rapha\"el Ponge}

\address{Department of Mathematics, University of Toronto, Canada.}
\email{ponge@math.toronto.edu}
\thanks{Research partially supported by NSF grant DMS 0409005 and JSPS Fellowship PE06016}


\begin{abstract}
    In this paper we produce several new invariants for CR and contact manifolds by looking at the noncommutative residue traces  
    of various geometric \psivdo\ projections. In the CR setting these operators arise from the $\dbarb$ complex and include the Szeg\"o projections 
    acting on forms. 
    In the contact setting they stem from  the generalized Szeg\"o projections at arbitrary integer levels of 
    Epstein-Melrose and from the contact complex of Rumin. In particular, we recover and extend recent results of 
    Hirachi and Boutet de Monvel and we answer a question of Fefferman.
\end{abstract} 
 
\maketitle 

\section{Introduction}
Motivated by Fefferman's program in CR geometry~\cite{Fe:PITCA}, Hirachi~\cite{Hi:LSSKGISPD} recently proved that the integral of the coefficient of the 
logarithmic singularity of the Szeg\"o kernel on the boundary of a strictly pseudoconvex domain in $\C^{n+1}$ gives rise to a CR invariant. 
This was subsequently extended to the contact setting by Boutet de Monvel~\cite{BdM:LTTP}
in terms of the Szeg\"o projections of~\cite{BG:STTO}. As later shown by Boutet de Monvel~\cite{BdM:RIMS05} these invariants vanish, but it
it was also  
asked by Fefferman~\cite{Fe:PC} whether there exist other examples of geometric operators such that the logarithmic singularities of their kernels give rise to CR 
or contact invariants.

The aim of this paper is to answer Fefferman's question by exhibiting various geometric projections on CR and contact manifolds such that the 
logarithmic singularities of their kernels give rise to invariants of the corresponding geometric structures. Furthermore, the framework that we used 
makes it possible to compute of the logarithmic singularities and the corresponding invariants by using techniques borrowed 
from index theory and Connes' noncommutative geometry. 

The Szeg\"o projection and its generalizations in~\cite{BG:STTO} are \psivdos\ in the sense of the Heisenberg calculus of~\cite{BG:CHM} and~\cite{Ta:NCMA}. 
Moreover, it has been shown by the 
author~(\cite{Po:CRAS1}, \cite{Po:GAFA1}) that the integral of the logarithmic singularity of the kernel of a \psivdo\ gives rise to a 
noncommutative residue trace for the Heisenberg calculus. Our invariants then appear as noncommutative residues of geometric \psivdo\ projections on CR and 
contact manifolds. These projections can be classified into three families of operators. 

%

The first family arises from the $\dbarb$-complex on CR manifolds. Namely, under $Y(q)$-type conditions the Szeg\"o projection on forms and
the orthogonal projections onto the kernels of the operators $\dbarb$ and $\dbarb^{*}$ are \psivdos\ (see, e.g.,~\cite{BG:CHM}). We then show 
that their noncommutative residues are all CR diffeomorphism invariants (Theorem~\ref{CR:Thm1}). In particular, in the strictly pseudoconvex case this 
allows us to recover Hirachi's result. The result further extends to include the projections 
associated to the $\dbarb$-complex with coefficients in a CR holomorphic vector bundle (Theorem~\ref{CR:Thm2}). 
In addition, we show that these invariants are not 
affected by deformations of the CR structure (see Propositions~\ref{prop:CR.deformation} and~\ref{prop:CR.deformation-spc}). 


The Szeg\"o projections of~\cite{BG:STTO} on a contact manifold $M$ have been further generalized by Epstein-Melrose~\cite{EM:HAITH} 
to arbitrary integer level and in such way to act on the sections of an arbitrary vector bundle $\cE$ over $M$. 
These operators are \psivdos\ and we show that the value of the noncommutative residue of a generalized Szeg\"o projection at a given 
integer level $k$ is independent of the 
choice of the operator and is an invariant of the Heisenberg diffeomorphism class of $M$ and of the $K$-theory class of $\cE$  (Theorem~\ref{thm:IGSP.main}).  
As a consequence this residue is  independent of the choice of the contact form and 
is invariant under deformations of the contact structure. Moreover,  
when $k=0$ and $\cE$ is the trivial line bundle this allows us to recover Boutet de Monvel's result. 

%
%
%

The last family of examples stems from the contact complex of Rumin~\cite{Ru:FDVC}. The latter is a complex of horizontal differential forms on 
a contact manifold which is hypoelliptic in every degree. The orthogonal projections onto the kernels of the 
differentials of this complex are \psivdos\ and we show that their noncommutative residues are Heisenberg diffeomorphism invariants and are invariant under deformation 
of the contact structure (Theorem~\ref{Contact:Thm}). 

The proofs for the examples arising from the $\dbarb$-complexe and the contact complex use simpler arguments than those 
of~\cite{Hi:LSSKGISPD} and~\cite{BdM:LTTP}, as the results follow from 
the observation that two \psivdo\ projections with same range or same kernel have same noncommutative residue (Lemma~\ref{prop:invariance-range-kernel}).
The proof for the examples coming from generalized Szeg\"o projections  partly relies on the fact that two \psivdo\ projections with homotopic 
principal symbols have same noncommutative residue (Lemma~\ref{lem:NCRPP.homotopy-invariance}). This generalizes the homotopy arguments 
of~\cite{Hi:LSSKGISPD} and~\cite{BdM:LTTP}. 

Next, the computation of these invariants is rather difficult. They appear as the integrals of local noncommutative residue densities, for which we have 
explicit formulas in terms of Heisenberg symbols. However, the number of terms to compute increases dramatically with the dimension, so there is no hope to 
get explicit geometric formulas without further tools to organize the computation. 

Furthermore, computating of the local noncommutative residue densities is even more important than the actual computating of the corresponding 
invariant, before the former implies the latter and could further provide us with some geometric 
information about the logarithmic singularities of the kernels of the corresponding \psivdo\ projections. 
At least in the case of the Szeg\"o kernel this would be of great interest in view of Fefferman's program. Therefore, even if the invariant may 
vanish it is interesting to compute the corresponding noncommutative residue densities.

In this paper we also allude to some new possible approaches for computing these densities and the corresponding invariants. 

A first approach that we suggest is to make use of Getzler's rescaling techniques in the setting of the Heisenberg calculus. This comes in naturally 
with the framework of the paper. It is believed that this could allow us to compute local densities associated to generalized Szeg\"o kernels, at least on 
strictly pseudoconvex CR manifolds (see Subsection~\ref{subsec:CR.computation} and Remark~\ref{rem:Szego.computation}).

%

Another approach suggested in Appendix is to make use of global $K$-theoretic techniques similar to those involved the $K$-theoretic proof of 
the index theorem of Atiyah-Singer~\cite{AS:IEO1}. To this end  we give 
a $K$-theoretic interpretation of the noncommutative residue of a \psivdo\ projection on a \emph{general} Heisenberg manifold 
$(M,H)$. More precisely, if we let $K_{0}(S_{0}(\fg^{*}M))$ denote the first $K$-group of the (noncommutative) algebra of zero'th order symbols, then 
the noncommutative residue of a \psivdo\ projection gives rise to an additive map $\rho_{R}:K_{0}(S_{0}(\fg^{*}M))\rightarrow \R$ (see 
Proposition~\ref{prop:NCR.rho-R}).

%

Notice that due to the noncommutativity of $S_{0}(\fg^{*}M)$ we really have to rely on the $K$-theory of algebras rather on that of spaces. 
Therefore, computing the map $\rho_{R}$ would definitely involve using tools from Connes' noncommutative geometry. As we also explain in Appendix 
two opposite interesting phenomena may occur:\smallskip 

(i)  The map $\rho_{R}$ is nontrivial and is computable in topological terms;\smallskip 

(ii) The map $\rho_{R}$ vanishes identically.\smallskip

Proving (i) could allow us to compute the invariants when we cannot use Getzler's rescaling techniques and this could be especially relevant for 
dealing with the invariants from the contact complex and with the CR invariants on CR manifolds with degenerate Levi form. However, the occurence of 
(ii) won't be too disappointing, because it would allow us to define the eta invariant of hypoelliptic selfadjoint \psivdos, which should be useful 
for dealing with index problems on complex manifolds with boundaries and on the asymptotically complex hyperbolic (ACH) manifolds. 

Finally, the arguments used in this paper are fairly general and should hold in many other settings as well. In particular, it would be interesting 
to extend them to the setting of complex manifolds with boundary and ACH manifolds. In particular, it would be of special interest to get an analogue in this context of 
Hirachi's invariant defined in terms of the Bergman projection.

The rest of the paper is organized as follows. In Section~\ref{sec.Heisenberg} we recall the main facts about the Heisenberg calculus and the noncommutative 
residue for this calculus. In Section~\ref{sec.NCRPP} we prove general results about noncommutative residues of \psivdo\ projections. 
In Section~\ref{sec:CR} we deal with the invariants from the $\dbarb$-complex on a CR manifold. 
Section~\ref{sec.Toeplitz} is devoted to the noncommutative residues of generalized Szeg\"o projections on a contact manifold. In Section~\ref{sec.contact-complex} 
we deal with the invariants arising from the contact complex. Finally, in Appendix we give a $K$-theoretic interpretation of the noncommutative 
residue of a \psivdo\ projection. 


\section{Heisenberg calculus and noncommutative residue}\label{sec.Heisenberg}
In this section we recall the main facts about the Heisenberg calculus and the noncommutative residue trace for this calculus. We also explain how 
the invariants of Hirachi and Boutet de Monvel can be interpreted as noncommutative residues. 
\subsection{Heisenberg manifolds}
A Heisenberg manifold is a pair $(M,H)$ consisting of a manifold $M^{d+1}$ together with a distinguished hyperplane bundle $H 
\subset TM$. This definition covers many examples: Heisenberg group, CR manifolds, contact manifolds, (codimension 1) foliations and the 
confoliations of~\cite{ET:C}. In addition, given another Heisenberg manifold $(M',H')$ we say that a diffeomorphism 
$\phi:M\rightarrow M'$ is a Heisenberg  diffeomorphism when $\phi_{*}H=H'$. 

The terminology Heisenberg manifold stems from the fact that the relevant tangent structure in this setting is that of a bundle $GM$ of graded nilpotent Lie 
groups (see~\cite{BG:CHM}, \cite{Be:TSSRG}, \cite{EMM:RLSPD}, \cite{FS:EDdbarbCAHG},  \cite{Gr:CCSSW}, \cite{Po:Pacific1}, \cite{Ro:INA}). 
This tangent Lie group bundle can be described as follows. 

First, we can define an intrinsic Levi form as the 2-form $\cL:H\times H\rightarrow TM/H$ such that, for any point $a 
\in M$ and any sections $X$ and $Y$ of $H$ near $a$, we have 
\begin{equation}
    \cL_{a}(X(a),Y(a))=[X,Y](a) \qquad \bmod H_{a}.
     \label{eq:NCRP.Levi-form}
\end{equation}
In other words the class of $[X,Y](a)$ modulo $H_{a}$ depends only on $X(a)$ and $Y(a)$, not on the germs of $X$ and $Y$ near $a$ (see~\cite{Po:Pacific1}). 

We define the tangent Lie algebra bundle $\fg M$ as the graded Lie algebra bundle consisting of $(TM/H)\oplus H$ together with the 
fields of Lie bracket and dilations such that, for sections $X_{0}$, $Y_{0}$ of $TM/H$ and $X'$, $Y'$ 
of $H$ and for $t\in \R$, we have 
\begin{equation}
    [X_{0}+X',Y_{0}+Y']=\cL(X',Y'), \qquad  t.(X_{0}+X')=t^{2}X_{0}+t X' .
    \label{eq:NCRP.Heisenberg-dilations} 
\end{equation}

Each fiber $\fg_{a}M$ is a two-step nilpotent Lie algebra, so by requiring the exponential map to be the identity 
the associated tangent Lie group bundle $GM$ appears as $(TM/H)\oplus H$ together with the grading above and the product law such that, 
for sections $X_{0}$, $Y_{0}$ of $TM/H$ and $X'$, $Y'$ of $H$, we have 
 \begin{equation}
     (X_{0}+X').(Y_{0}+Y')=X_{0}+Y_{0}+\frac{1}{2}\cL(X',Y')+X'+Y'.
 \end{equation}
%
%
%
%

 Moreover, if $\phi$ is a Heisenberg diffeomorphism from $(M,H)$ onto a Heisenberg manifold $(M',H')$ then, as $\phi_{*}H=H'$ we get linear isomorphisms 
 from $TM/H$ onto $TM'/H'$ and from $H$ onto $H'$ which together give rise to a linear isomorphism 
 $\phi_{H}':(TM/H)\oplus H\rightarrow (TM'/H')\oplus H'$. In fact $\phi_{H}'$ is a graded Lie group isomorphism from $GM$ onto $GM'$ 
 (see~\cite{Po:Pacific1}). 
 
 \subsection{Heisenberg calculus}  
 The Heisenberg calculus is the relevant pseudodifferential calculus to study hypoelliptic 
 operators on Heisenberg manifolds. It was independently introduced by Beals-Greiner~\cite{BG:CHM} and Taylor~\cite{Ta:NCMA} (see 
 also~\cite{BdM:HODCRPDO}, \cite{Dy:POHG}, \cite{Dy:APOHSC}, \cite{EM:HAITH}, \cite{FS:EDdbarbCAHG}, \cite{RS:HDONG}).
  
 The initial idea in the Heisenberg calculus, which is due to Stein, is to construct a 
 class of operators on a Heisenberg manifold $(M^{d+1},H)$, called \psivdos, which at any point $a \in M$ are modeled on homogeneous left-invariant 
 convolution operators on the tangent group $G_{a}M$. 
%
 
 Locally the \psidos\ can be described  as follows. Let $U \subset \Rd$ be a local chart together with a frame $X_{0},\ldots,X_{d}$ of $TU$ such that 
 $X_{1},\ldots,X_{d}$ span $H$. Such a chart is called a Heisenberg chart. Moreover, on $\Rd$ we consider the dilations, 
\begin{equation}
     t.\xi=(t^{2}\xi_{0},t\xi_{1},\ldots,t\xi), \qquad \xi\in \Rd, \quad t>0. 
\end{equation}
 
\begin{definition}1)  $S_{m}(\URd)$, $m\in\C$, is the space of functions 
    $p(x,\xi)$ in $C^{\infty}(U\times\Rdo)$ such that $p(x,t.\xi)=t^m p(x,\xi)$ for any $t>0$.\smallskip

2) $S^m(\URd)$,  $m\in\C$, consists of functions  $p\in C^{\infty}(\URd)$ with
an asymptotic expansion $ p \sim \sum_{j\geq 0} p_{m-j}$, $p_{k}\in S_{k}(\URd)$, in the sense that, for any integer $N$ and 
for any compact $K \subset U$, we have
\begin{equation}
    | \partial^\alpha_{x}\partial^\beta_{\xi}(p-\sum_{j<N}p_{m-j})(x,\xi)| \leq 
    C_{\alpha\beta NK}\|\xi\|^{\Re m-\brak\beta -N}, \qquad  x\in K, \quad \|\xi \| \geq 1,
    \label{eq:NCRP.asymptotic-expansion-symbols}
\end{equation}
where we have let $\brak\beta=2\beta_{0}+\beta_{1}+\ldots+\beta_{d}$ and 
 $\|\xi\|=(\xi_{0}^{2}+\xi_{1}^{4}+\ldots+\xi_{d}^{4})^{1/4}$.
\end{definition}

Next, for $j=0,\ldots,d$ let  $\sigma_{j}(x,\xi)$ denote the symbol (in the 
classical sense) of the vector field $\frac{1}{i}X_{j}$  and set  $\sigma=(\sigma_{0},\ldots,\sigma_{d})$. Then for $p \in S^{m}(\URd)$ we let $p(x,-iX)$ be the 
continuous linear operator from $C^{\infty}_{c}(U)$ to $C^{\infty}(U)$ such that 
    \begin{equation}
          p(x,-iX)f(x)= (2\pi)^{-(d+1)} \int e^{ix.\xi} p(x,\sigma(x,\xi))\hat{f}(\xi)d\xi,
    \qquad f\in C^{\infty}_{c}(U).
    \end{equation}


\begin{definition}
   $\pvdo^{m}(U)$, $m\in \C$, consists of operators $P:C^{\infty}_{c}(U)\rightarrow C^{\infty}(U)$ which are of the form
$P= p(x,-iX)+R$ for some $p$ in $S^{m}(\URd)$, called the symbol of $P$, and some smoothing operator $R$.
\end{definition}

 For any $a\in U$ there is exists a unique affine change of variable $\psi_{a}:\Rd \rightarrow \Rd$ such that $\psi_{a}(a)=0$ 
 and $(\psi_{a})_{*}X_{j}=\frac{\partial}{\partial x_{j}}$ at $x=0$ for $j=0,1,\ldots,d+1$. Then, a continuous operator 
 $P:C^{\infty}_{c}(U)\rightarrow C^{\infty}(U)$ is a \psivdo\ of order $m$ if, and only if, its kernel $k_{P}(x,y)$ has a behavior near the diagonal 
 of the form, 
 \begin{equation}
     k_{P}(x,y) \sim \!\!\sum_{j\geq -(m+d+2)} \!\! (a_{j}(x,\psi_{x}(y))- \!\! \sum_{\brak \alpha =j}\!\! c_{\alpha}(x)\psi_{x}(x)^{\alpha}\log 
     \|\psi_{x}(y)\|), 
      \label{eq:NCRP.kernel-characterization}
 \end{equation}
 with $c_{\alpha}\in C^{\infty}(U)$ and $a_{j}(x,y) \in C^{\infty}(\URdo)$ such that 
 $a_{j}(x,\lambda.y)=\lambda^{j}a_{j}(x,y)$ for any $\lambda>0$. Moreover, $a_{j}(x,y)$ and $c_{\alpha}(x)$, $\brak \alpha=j$, depend only on the symbol 
 of $P$ of degree $-(j+d+2)$.

 
 The class of \psivdos\ is invariant under changes of Heisenberg chart (see~\cite[Sect.~16]{BG:CHM}, \cite[Appendix A]{Po:MAMS1}), so we may 
 extend the definition of \psivdos\ to an arbitrary Heisenberg manifold $(M,H)$ and let them act on sections of a vector bundle $\cE$ over $M$. 
 We let $\pvdo^{m}(M,\cE)$ denote the class of \psivdos\ of order $m$ on $M$ acting on sections 
 of $\cE$. 
 
 From now on we let $(M^{d+1},H)$ be a compact Heisenberg manifold and we let $\fg^{*}M$ denote the (linear) dual of the Lie algebra bundle $\fg M$ of $GM$ 
 with canonical projection  $\text{pr}: M\rightarrow \fg^{*}M$. As shown 
 in~\cite{Po:MAMS1} (see also~\cite{EM:HAITH}) 
 the principal symbol of $P\in \pvdo^{m}(M,\cE)$ can be intrinsically defined as a symbol $\sigma_{m}(P)$ of the class below. 
 
\begin{definition}
  $S_{m}(\fg^{*}M)$, $m\in \C$, consists of sections $p\in C^{\infty}(\fg^{*}M\setminus 0, \End \textup{pr}^{*}\cE)$ which are homogeneous of 
    degree $m$ with respect to the dilations in~(\ref{eq:NCRP.Heisenberg-dilations}), i.e., we have 
   $p(x,\lambda.\xi)=\lambda^{m}p(x,\xi)$ for any $\lambda>0$. 
\end{definition}

 Next, for any $a \in M$ the convolution 
 on $G_{a}M$ gives rise under the (linear) Fourier transform to a bilinear product for homogeneous symbols, 
 \begin{equation}
     *^{a}: S_{m_{1}}(\fg^{*}_{a}M,\cE_{a})\times S_{m_{2}}(\fg^{*}_{a}M,\cE_{a}) \longrightarrow S_{m_{1}+m_{2}}(\fg^{*}_{a}M,\cE_{a}),
 \end{equation}
 This product depends smoothly on $a$ as much so to yield a product, 
 \begin{gather}
     *: S_{m_{1}}(\fg^{*}M,\cE)\times S_{m_{2}}(\fg^{*}M,\cE) \longrightarrow S_{m_{1}+m_{2}}(\fg^{*}M,\cE),\\
     p_{m_{1}}*p_{m_{2}}(a,\xi)=[p_{m_{1}}(a,.)*^{a}p_{m_{2}}(a,.)](\xi).
     \label{eq:NCRP.product-symbols}
 \end{gather}
 This provides us with the right composition for principal symbols, for we have
 \begin{equation}
     \sigma_{m_{1}+m_{2}}(P_{1}P_{2})=\sigma_{m_{1}}(P_{1})*\sigma_{m_{2}}(P_{2}) \qquad \forall P_{j}\in \pvdo^{m_{j}}(M,\cE).
 \end{equation}

 Notice that when $G_{a}M$ is not commutative, i.e., $\cL_{a}\neq 0$, the product $*^{a}$ is not anymore the pointwise product of symbols and, in 
 particular, is not commutative. 
 Consequently, unless when $H$ is integrable, the product for Heisenberg symbols is not commutative and, while local, it is not microlocal. 
 
 If $P\in \pvdo^{m}(M,\cE)$ then the transpose $P^{t}$ belongs to $\pvdo^{m}(M,\cE^{*})$, and if $\cE$ is further endowed with a Hermitian 
 metric then the adjoint $P^{*}$ belongs to $\pvdo^{\overline{m}}(M,\cE)$ (see~\cite{BG:CHM}). Moreover, as shown in~\cite[Sect.~3.2]{Po:MAMS1}, 
 the principal symbols of $P^{t}$ and $P^{*}$ are 
 \begin{equation}
     \sigma_{m}(P^{t})=\sigma_{m}(P)(x,-\xi)^{t} \quad \text{and} \quad \sigma_{\overline{m}}(P^{*})=\sigma_{m}(P)(x,\xi)^{*}.
      \label{eq:HC.symbols-transpose-adjoints}
 \end{equation}
 
 When the principal symbol of $P\in \pvdo^{m}(M,\cE)$ is invertible with respect to the product $*$, the symbolic calculus 
 of~\cite{BG:CHM} allows us to construct a parametrix for $P$ in $\pvdo^{-m}(M,\cE)$. In particular, although not elliptic, 
 $P$ is hypoelliptic with a controlled loss/gain of derivatives (see~\cite{BG:CHM}). 
 
 In general, it may be difficult to determine whether the principal symbol of a given \psivdo\ $P\in \pvdo^{m}(m,\cE)$ is invertible with respect to the product 
 $*$, but  this can be completely determined in terms of a representation theoretic criterion on each tangent 
 group $G_{a}M$, the so-called Rockland condition (see~\cite[Thm.~3.3.19]{Po:MAMS1}).  In particular, 
 if $\sigma_{m}(P)(a,.)$ is \emph{pointwise} invertible with respect to the product $*^{a}$ for any 
 $a \in M$ then  $\sigma_{m}(P)$ is \emph{globally} invertible with respect to $*$.

  \subsection{Noncommutative residue}
 Let $P:C^{\infty}(M,\cE)\rightarrow C^{\infty}(M,\cE)$ be a \psivdo\ of integer order $m$. Then it follows from~(\ref{eq:NCRP.kernel-characterization}) 
 that in a trivializing Heisenberg chart the kernel $k_{P}(x,y)$ of $P$ has a behavior near the diagonal of the form,  
      \begin{equation}
        k_{P}(x,y)=\!\!\sum_{-(m+d+2)\leq j\leq 1}\!\!a_{j}(x,-\psi_{x}(y)) - c_{P}(x)\log \|\psi_{x}(y)\| + 
        \op{O}(1),
        \label{eq:NCRHC.log-singularity}
    \end{equation}
    where $a_{j}(x,y)$ is homogeneous of degree $j$ in $y$ with respect to the dilations~(\ref{eq:NCRP.Heisenberg-dilations}) 
    and $c_{P}(x)$ is the smooth function given by 
    \begin{equation}
        c_{P}(x)=|\psi_{x}'|\int_{\|\xi\|=1}p_{-(d+2)}(x,\xi)d\xi,
         \label{eq:NCRHC.formula-cP}
    \end{equation}
    where $p_{-(d+2)}(x,\xi)$ is the homogeneous symbol of degree $-(d+2)$ of $P$. 
    
   Under the action of Heisenberg diffeomorphisms $c_{P}(x)$ behaves like a density (see~\cite[Prop.~3.11]{Po:GAFA1}). 
   Therefore, the coefficient $c_{P}(x)$ makes intrinsically sense on $M$ as a section of $|\Lambda|(M)\otimes \End \cE$, where $|\Lambda|(M)$ 
    is the bundle of densities on $M$. 
   
    We can now define a functional on $\pvdoz(M,\cE)=\cup_{m \in \Z}\pvdo^{m}(M,\cE)$ by letting 
    \begin{equation}
        \Res P= \int_{M} \tr_{\cE} c_{P}(x), \qquad P\in \pvdoz(M,\cE). 
    \end{equation}

    As shown in~\cite{Po:GAFA1} this functional is the analogue for the Heisenberg calculus of the noncommutative residue of 
Wodzicki~(\cite{Wo:LISA}, \cite{Wo:NCRF}) and Guillemin~\cite{Gu:NPWF}, since it also arises as the residual trace on integer order 
\psivdos\ induced by the analytic continuation of the usual trace to \psivdos\ of non-integer orders. 

  \begin{proposition}[\cite{Po:CRAS1}, \cite{Po:GAFA1}]\label{prop:NCR}
     1) $\Res$ is a trace on the algebra $\pvdoz(M,\cE)$ which vanishes on differential operators and on \psivdos\ of integer order~$\leq -(d+3)$, 
     In fact, when $M$ is connected this is the unique trace up to constant multiple.\smallskip
     
     2) For any $P \in \pvdoz(M,\cE)$ we have $\Res P^{t}=\overline{\Res P^{*}}=\Res P$.\smallskip
     
     3) Let $\phi$ be a Heisenberg diffeomorphism from $(M,H)$ onto a Heisenberg manifold $(M',H')$. Then for any $P \in \pvdoz(M,\cE)$ we have $\Res 
     \phi_{*}P=\Res P$. 
 \end{proposition}

 \subsection{Logarithmic singularity of Szeg\"o kernels}
%
Let $S$ be a  Szeg\"o projection on a contact manifold $M^{2n+1}$ as in~\cite{BG:STTO}. This is a FIO with complex phase $q(x,y)$ and 
near the diagonal the kernel of $S$ a behavior of the form 
\begin{equation}
    k_{S}(x,y)\sim \sum_{-(n+1)\leq j\leq -1}\alpha_{j}(x,y)q(x,y)^{j} +\sum_{j\geq 0}\beta_{j}(x,y)q(x,y)^{j}\log q(x,y),
    \label{eq:NRHC.singularity-Szego}
\end{equation}
where $\alpha_{j}(x,y)$ and $\beta_{j}(x,y)$ are smooth functions defined near the diagonal. 

The coefficient $\beta_{0}(x,x)$ of the logarithmic singularity makes sense globally as a density on $M$ and so we can define 
\begin{equation}
    L(S)=\int_{M}\beta_{0}(x,x). 
\end{equation}
This is this object which is shown to give rise to a global invariant in~\cite{Hi:LSSKGISPD} and~\cite{BdM:LTTP}.

In fact, the phase $q(x,y)$ vanishes on the diagonal and is such that 
$id_{x}q=-id_{y}q$ is a nonzero annihilator of $H$ on the diagonal and $\Re q(x,y)\gtrsim |x-y|^{2}$  near the diagonal. Therefore, the Taylor expansion of 
$q(x,y)$ near $y=x$ is of the form, 
\begin{equation}
    q(x,y) \sim \sum_{\brak \alpha\geq 2} a_{\alpha}(x)\psi_{x}(y)^{\alpha}, \qquad 
     \label{eq:NRHC.Taylor-complex-phase}
\end{equation}
where $q_{2}(x,y):=\sum_{\brak \alpha\geq 2} a_{\alpha}(x)\psi_{x}(y)^{\alpha}$ is nonzero for $y\neq 0$ and $y$ close enough to $x$. 

In fact, plugging~(\ref{eq:NRHC.Taylor-complex-phase}) into~(\ref{eq:NRHC.singularity-Szego}) shows that the kernel of $S$ has near the diagonal the
singularity of the form~(\ref{eq:NCRP.kernel-characterization}) with $m=0$ and so $S$ is a zero'th order \psivdo. Moreover, 
as near the diagonal we have $\log q(x,y)=2 \log \|\psi_{x}(y)\|+\op{O}(1)$,   we see that
$c_{S}(x)=-2\beta_{0}(x,x)$. Thus, 
\begin{equation}
    L(S)=-\frac{1}{2} \Res S.
\end{equation}
This shows that the invariants considered by Hirachi and Boutet de Monvel can be interpreted as noncommutative residues. 

\begin{remark}\label{rem:NCR.FIO}
    Guillemin~\cite{Gu:RTCAFIO} has defined noncommutative residue traces for some algebras of FIO's, including the algebra of T\"oplitz operators 
    on a contact manifold. The latter is an ideal of the algebra of \psivdos\ (see, e.g.,~\cite{EM:HAITH}) and one can check that in this context 
    Guillemin's trace is equal to $-\frac{1}{2}\Res $ on T\"oplitz operators. In particular, we see that $L(S)$ agrees with the noncommutative residue trace in 
    Guillemin's sense of $S$. 
\end{remark}

\section{Noncommutative residues of \psivdo\ projections}\label{sec.NCRPP}
This section  we gather several general lemmas about noncommutative residues of a \psivdo\ projections. 

Throughout all the section we let $(M^{d+1},H)$ be a compact Heisenberg manifold equipped with a smooth 
density~$>0$ and let $\cE$ be a Hermitian vector bundle. 

%
%
%
%
\begin{lemma}\label{lem:similarity-lemma}
    Let $\Pi\in \pvdo^{0}(M,\cE)$ be a \psivdo\ projection. Then the orthogonal projection $\Pi_{0}$ onto 
    its range is a zero'th order \psivdo\ and we have $\Res \Pi_{0}=\Res \Pi$.
\end{lemma}
\begin{proof}
    It is well known that any projection on $L^{2}(M,\cE)$ is similar to the orthogonal projection onto its range (see, e.g.,~\cite[Prop.~4.6.2]{Bl:KTOA}). 
    Indeed, the operator $B=1+(\Pi-\Pi^{*})(\Pi^{*}-\Pi)$ is invertible on $L^{2}(M,\cE)$ and we have $  \Pi_{0}=\Pi\Pi^{*}B^{-1}$.
    Moreover, if we let $A=\Pi_{0}+(1-\Pi)(1-\Pi_{0})$ then $A$ is invertible with inverse $A^{-1}=\Pi+(1-\Pi_{0})(1-\Pi)$ and we have $ 
    \Pi_{0}=A^{-1}\Pi A$.
%
    
    Observe that $B$ is a zero'th order \psivdo. Let $Q= \Pi-\Pi^{*}$ and for $a\in M$  let $B^{a}$ and 
    $Q^{a}$ be the respective model operators at $a$ of $B$ as defined in~\cite[Sect.~3.2]{Po:MAMS1}. Recall that the latter are bounded
    left invariant convolution operators on $L^{2}(G_{a}M,\cE_{a})$. 
    Since $B=1+QQ^{*}$ it follows from~\cite[Props.~3.2.9, 3.2.12]{Po:MAMS1} that, either 
    $B^{a}=1+Q^{a}(Q^{a})^{*}$ if $Q$ has order 0, or $B^{a}=1$ if $Q$ otherwise.
    
    In any case $B^{a}$ is an invertible bounded operator on $L^{2}(G_{a}M,\cE)$, so it follows from~\cite[Thm.~3.3.10]{Po:MAMS1} 
    that the principal symbol of $B$ is invertible. Therefore, there exist $C\in \pvdo^{0}(M,\cE)$ and smoothing operators $R_{1}$ and $R_{2}$ such that
    $CB=1-R_{1}$ and $BC=1-R_{2}$.
    From this we get $C=B^{-1}-R_{1}B^{-1}=B^{-1}-B^{-1}R_{2}$. This implies that $B^{-1}$ induces a continuous endomorphism of 
    $C^{\infty}(M,\cE)$ and agrees with $C$ modulo a smoothing operator, hence is a zero'th order \psivdo. 
    
    Since $B^{-1}$ is a zero'th order \psivdo\ we deduce that $\Pi_{0}=\Pi\Pi^{*}B^{-1}$ is a zero'th order \psivdo\ as well. This implies that $A$ 
    and $A^{-1}$ are also \psivdos, so as $\Res$ is a trace we get $\Res \Pi_{0}=\Res A^{-1}\Pi A=\Res \Pi$. The lemma is thus proved.
\end{proof}

As a consequence of this lemma we will obtain:

\begin{lemma}\label{prop:invariance-range-kernel}
   For $j=1,2$ let $\Pi_{j}\in \pvdo^{0}(M,\cE)$ be a \psivdo\ projection. If $\Pi_{1}$ and $\Pi_{2}$ have same range or have same kernel then $\Res 
   \Pi_{1}=\Res \Pi_{2}$. 
\end{lemma}
\begin{proof}
    If $\Pi_{1}$ and $\Pi_{2}$ have same range then by Lemma~\ref{lem:similarity-lemma} their noncommutative residues agree, since their are both 
    equal to that of the orthogonal projection onto their common range. 
    
    If $\Pi_{1}$ and $\Pi_{2}$ have same kernel then $1-\Pi_{1}$ and $1-\Pi_{2}$ have same range, so we have $\Res (1-\Pi_{2})=\Res 
    (1-\Pi_{2})$. As $\Res (1-\Pi_{j})=-\Res \Pi_{j}$, $j=1,2$, it follows that we have $\Res \Pi_{2}=\Res \Pi_{1}$.
\end{proof}

Another consequence of Lemma~\ref{lem:similarity-lemma}  is the following.

\begin{lemma}\label{NCRP.real-valued}
    Let $\Pi_{0}\in \pvdo^{0}(M,\cE)$ be a \psivdo\ projection. Then $\Res \Pi$ is a real number and we have $\Res \Pi^{*}=\Res \Pi^{t}=\Res \Pi$. 
\end{lemma}
\begin{proof}
   By Proposition~\ref{prop:NCR} we have $\overline{\Res \Pi^{*}}=\Res \Pi^{t}=\Res \Pi$, so we only have to check that $\Res \Pi$ is in $\R$.
   Let $\Pi_{0}$ be the orthogonal projection onto the range of $\Pi$. As $\Pi_{0}$ is a selfadjoint \psivdo\ projection we have $\overline{\Res 
   \Pi_{0}}=\Res \Pi_{0}^{*}=\Res \Pi_{0}$, so that $\Res \Pi_{0}$ is a real number. Since by Lemma~\ref{lem:similarity-lemma} the latter agrees with 
   $\Res \Pi$, we see that $\Res \Pi$ is in $\R$ as well. 
\end{proof}

Next, we define $C^{1}$-paths of \psivdos\ as follows. For an open $V\subset \Rd$ we endow $S^{m}(V\times \Rd)$, $m \in \C$, with the Fr\'echet space 
topology induced by the topology of $C^{\infty}(V\times \Rd)$ and the sharpest constants 
in~(\ref{eq:NCRP.asymptotic-expansion-symbols}). We then let $S^{m}(V\times \Rd)_{t}$ denote the space of $C^{1}$-paths from $I:=[0,1]$ to $S^{m}(V\times \Rd)$. 

Similarly 
we endow $S_{m}(\fg^{*}M, \cE)$ with the Fr\'echet space topology inherited from that of $C^{\infty}(\fg^{*}M\setminus 0, \End \op{pr}^{*}\cE)$ and we 
let $S_{m}(\fg^{*}M, \cE)_{t}$ denote the space of $C^{1}$-paths from $I$ to $S_{m}(\fg^{*}M, \cE)$.

\begin{definition}
    $\pvdo^{m}(M,\cE)_{t}$ is the space of paths $(P_{t})_{t\in I}\subset \pvdo^{m}(M,\cE)$ which are $C^{1}$ in the sense 
    that:\smallskip 
    
    (i) The kernel of $P_{t}$ is given outside the diagonal by a $C^{1}$-path of smooth kernels;\smallskip

(ii) For any Heisenberg chart $\kappa:U \rightarrow V\subset \Rd$ with a $H$-frame $X_{0},\ldots,X_{d}$ and any trivialization 
$\tau:\cE_{|_{U}}\rightarrow U\times \C^{r}$ we can write 
\begin{equation}
    \kappa_{*}\tau_{*}(P_{t|_{U}})=p_{t}(x,-iX)+R_{t},
\end{equation}
for some $C^{1}$-path  $p_{t}\in M_{r}(S^{m}(U\times \Rd)_{t})$ and some $C^{1}$-path $R_{t}$ of smoothing operators, i.e., 
$R_{t}$ is given by a $C^{1}$-path of 
smooth kernels.
\end{definition}

We gather the main properties of $C^{1}$-paths of \psivdos\ in the following.

\begin{lemma}[{\cite[Chap.~4]{Po:MAMS1}}]\label{lem:NCRPP-C1-path}
    1) If $P_{t}\in \pvdo^{m}(M,\cE)_{t}$ 
then $\sigma_{m}(P_{t})$  is  a $C^{1}$-path with values in $S_{m}(\fg^{*}M, \cE)$ and, in fact, 
in a local trivializing chart all the homogeneous components of the symbol of $P_{t}$ yield $C^{1}$-paths of homogeneous symbols.\smallskip 

2)    If $P_{j,t}\in \pvdo^{m_{j}}(M,\cE)_{t}$, $j=1,2$, then $P_{1,t}P_{2,t}\in \pvdo^{m_{1}+m_{2}}(M,\cE)_{t}$.\smallskip 

3) Let $\phi$ be a Heisenberg diffeomorphism from a Heisenberg manifold $(M',H')$ onto $(M,H)$. Then for any $P_{t}\in \pvdo^{m}(M,\cE)_{t}$ the 
path $\phi^{*}P_{t}$ is in $\pvdo^{m}(M',\phi^{*}\cE)_{t}$.
\end{lemma}

\begin{remark}
   In~\cite{Po:MAMS1} the proofs are actually carried out for holomorphic families of \psidos, but they remain valid \emph{mutatis mutandis} for 
   $C^{1}$-paths of \psivdos. 
\end{remark}

%
%
%

Bearing this in mind we have:

\begin{lemma}\label{lem:NCRPP.homotopy-invariance}
     Let $\Pi_{0}$ and $\Pi_{1}$ be  projections in $\pvdo^{0}(M,\cE)$ such that  their principal symbols can be joined to each other by means of 
     a $C^{1}$ path of idempotents in 
 $S_{0}(\fg^{*}M,\cE)$. Then we have $\Res \Pi_{0}=\Res \Pi_{1}$.
\end{lemma}
\begin{proof}
 For $j=0,1$ let $F_{j}=2\Pi_{j}-1$. Then $F_{j}^{2}=1$ and the principal symbol of $F_{0}$ can be connected to that of 
 $F_{1}$ by means of a $C^{1}$-path $f_{0,t}\in S_{0}(\fg^{*}M,\cE)_{t}$ such that $f_{0,t}*f_{0,t}=1$  $\forall t\in[0,1]$. 
  We can 
 construct a path $G_{t}\in \pvdo^{0}(M,\cE)_{t}$ so that we have $G_{t}^{2}=1 \bmod \pvdo^{-(d+3)}(M,\cE)_{t}$ and $G_{j}=F_{j} \bmod \pvdo^{-(d+3)}(M,\cE)$ 
      for $j=0,1$ as follows. 
 
      Let $(\varphi_{k})_{k\geq 0} \subset C^{\infty}(M)$ be a partition of the unity subordinated to an open covering $(U_{k})_{k\geq 0}$ of domains of 
     Heisenberg charts $\kappa_{k}:U_{k}\rightarrow V_{k}$ with $H$-frame $X_{0}^{(k)},\ldots,X_{d}^{(k)}$ and 
     over which there are trivializations $\tau_{k}:\cE_{|_{U_{k}}}\rightarrow U_{k}\times 
     \C^{r}$.  Let $q_{0,t}^{(k)}(x,\xi)=(1-\chi(\xi))(\kappa_{\alpha*}\tau_{\alpha*}f_{t})(x,\xi)$
     where $\chi\in C^{\infty}_{c}(\Rd)$ is such that $\chi(\xi)=1$ near the origin. Then $q_{0,t}^{(k)}$ is a $C^{1}$-path with values in $S^{0}(U_{k}\times 
     \Rd,\C^{r})$ and we get a path $Q_{t}\in \pvdo^{0}(M,\cE)_{t}$ with principal symbol $f_{0,t}$ by letting 
     \begin{equation}
         Q_{t}=\sum_{k\geq 0}\varphi_{k}(\tau_{k}^{*}\kappa_{k}^{*}q_{0,t}(x,-iX^{(k)}))\psi_{k},
     \end{equation}
     where $\psi_{k}\in C^{\infty}_{c}(U_{k})$ is such that $\psi_{k}=1$ near $\op{supp}\varphi_{k}$. 
     
     Since $Q_{t}$ has principal symbol $f_{0,t}$ we see that $F_{j}-Q_{j}\in \pvdo^{-1}(M,\cE)$ for $j=0,1$. Consider now the $C^{1}$-path,
    $P_{t}=Q_{t}+(1-t)(F_{0}-Q_{0})+t(F_{1}-Q_{1})$. 
     Then $P_{t}$ has principal symbol $f_{0,t}$ for every $t\in [0,1]$ and for $j=0,1$ we have $P_{j}=F_{j}$. 
%
%
  
     Next, since 
      $f_{0,t}*f_{0,t}=1$ we can write $P_{t}^{2}=1-R_{t}$ with $R_{t}\in \pvdo^{-1}(M,\cE)_{t}$. Let $\sum_{k\geq 0}a_{k}z^{k}$ be the 
      Taylor series at $z=0$ of $(1-z)^{-\frac{1}{2}}$. In particular, for any integer $N$ we can write  
      $(1-z)(\sum_{0\leq k\leq N} a_{k}z^{k})^{2}=1+P_{N}(z)$, where $P_{N}(z)$ is a polynomial of the form 
      $P_{N}(z)=\sum_{N+1\leq k+l\leq 2N+1}b_{N,k}z^{k}$. 
      
      Let $G_{t}=P_{t}\sum_{0\leq k\leq d+2}a_{k}R_{t}^{k}$. Then $G_{t}\in \pvdo^{0}(M,\cE)_{t}$ and as $R_{t}=1-P_{t}^{2}$ commutes with $P_{t}$ 
      we get
  \begin{equation}
    G_{t}^{2}=P_{t}^{2}(\sum_{0\leq k\leq d+2} a_{k}R_{t}^{k})^{2}=(1-R_{t})(\sum_{0\leq k\leq d+2} a_{k}R_{t}^{k})^{2}=1+P_{N}(R_{t}). 
 \end{equation}
 Since $P_{d+2}(R_{t})$ is in $ \pvdo^{-(d+3)}(M,\cE)_{t}$
 this shows that $G_{t}^{2}=1 \bmod 
 \pvdo^{-(d+3)}(M,\cE)_{t}$.  Moreover, as for $j=0,1$ we have $R_{j}=1-F_{j}^{2}=0$ we see that $G_{j}=P_{j}=F_{j}$. 
    
 Now, the equality $G_{t}^{2}=1 \bmod \pvdo^{-(d+3)}(M,\cE)_{t}$ implies that $\dot{G}_{t}G_{t}+G_{t}\dot{G}_{t}=0$ modulo   
 $\pvdo^{-(d+3)}(M,\cE)_{t}$. Therefore $-G_{t}\dot{G}_{t}G_{t}=G_{t}^{2}\dot{G}_{t}=\dot{G}_{t} \bmod \pvdo^{-(d+3)}(M,\cE)_{t}$, and so we get
 $  \dot{G}_{t}=\frac{1}{2}(G_{t}^{2}\dot{G}_{t}-G_{t}\dot{G}_{t}G_{t})= \frac{1}{2}[G_{t},G_{t}\dot{G}_{t}] \bmod \pvdo^{-(d+3)}(M,\cE)_{t}$.

 
 On the other hand, it follows from~(\ref{eq:NCRHC.formula-cP}) and Lemma~\ref{lem:NCRPP-C1-path} 
 that $\Res $ commutes with the differentiation of $C^{1}$-paths. Since $\Res$  is a trace and 
 vanishes on $\pvdo^{-(d+3)}(M,\cE)$,  we obtain $  \frac{d}{dt}\Res G_{t}=\Res \dot{G}_{t}=\frac{1}{2}\Res [G_{t},G_{t}\dot{G}_{t}]=0$.
%
%
 Hence we have $\Res G_{0}=\Res G_{1}$. As $\Res G_{j}=\Res F_{j}=\Res (2\Pi_{j}-1)=2\Res \Pi_{j}$  it follows that $\Res \Pi_{0}=\Res \Pi_{1}$ as desired. 
 \end{proof}

\section{Invariants from the $\overline{\partial}_{b}$-complex}\label{sec:CR}
Throughout all this section we let $M^{2n+1}$ be a compact orientable CR manifold with CR bundle $T_{1,0}\subset T_{\C}M$, 
so that $H=\Re (T_{1,0}\oplus T_{0,1})\subset TM$ is a hyperplane bundle of $TM$ admitting an (integrable) complex structure. 

\subsection{Construction of the CR invariants}
Since $M$ is orientable and $H$ is orientable by means of its complex structure, there exists a global non-zero real 1-form $\theta$ annihilating $H$.  
Associated to $\theta$ is the Hermitian Levi form,  
 \begin{equation}
    L_{\theta}(Z,W)=-id\theta(Z,\overline{W})=i\theta([Z,\overline{W}]), \qquad Z,W\in C^{\infty}(M,T_{1,0}). 
      \label{eq:CR.Levi-form}
 \end{equation}
We then say that $M$ is strictly pseudoconvex (resp.~$\kappa$-strictly pseudoconvex) if for some choice of $\theta$ the Levi form is everywhere 
positive definite (resp.~has everywhere $\kappa$ negative eigenvalues and $n-\kappa$ positive eigenvalues). 

Let $\cN$ be a supplement of $H$ in $TM$. This is an orientable line bundle which gives rise to the splitting, 
\begin{equation}
    T_{\C}M=T_{1,0}\oplus T_{0,1}\oplus (\cN\otimes \C).
    \label{eq:CR-decomposition}
\end{equation}
For $p,q=0,\ldots,n$ let $\Lambda^{p,q}=(\Lambda^{1,0})^{p}\wedge (\Lambda^{0,1})^{q}$ be the bundle of $(p,q)$-covectors, where 
$\Lambda^{1,0}$ (resp.~$\Lambda^{0,1}$) denotes the annihilator 
in $T^{*}_{\C}M$ of $T_{0,1}\oplus (\cN\otimes \C)$ (resp.~of $T_{1,0}\oplus (\cN\otimes \C)$). Then we have the splitting, 
\begin{equation}
    \Lambda^{*}T_{\C}^{*}M=(\bigoplus_{p,q=0}^{n}\Lambda^{p,q})\oplus \theta\wedge  \Lambda^{*}T_{\C}^{*}M.
     \label{eq:CR-Lambda-pq-decomposition}
\end{equation}
Notice that this decomposition does not depend on the choice of $\theta$, but it does depend on that of $\cN$. 

The complex $\dbarb:C^{\infty}(M, \Lambda^{0,*})\rightarrow C^{\infty}(M,\Lambda^{0,*+1})$ of 
Kohn-Rossi~(\cite{KR:EHFBCM},~\cite{Ko:BCM}) is defined as follows. For any $\eta \in C^{\infty}(M, \Lambda^{0,q})$ we can uniquely decompose $d\eta$ as 
\begin{equation}
    d\eta =\dbarbq \eta + \partial_{b;q}\eta + \theta \wedge \cL_{X_{0}}\eta,
     \label{eq:CR.dbarb}
\end{equation}
where $\dbarbq \eta $  and $\partial_{b;q}\eta$ are sections of $\Lambda^{0,q+1}$ and $\Lambda^{1,q}$ respectively and $X_{0}$ is the section of $\cN$ 
such that $\theta(X_{0})=1$. Thanks to the integrability of $T_{1,0}$ we have $\overline{\partial}_{b;q+1}\circ \dbarbq=0$, so we get a cochain 
complex. 

The $\dbarb$-complex depends only on the CR structure of $M$ and on the choice of $\cN$. The dependence on the latter can be determined as follows. 
Let $\cN'$ be another supplement of $H$ and let us assign the superscipt $'$ to objects defined 
using $\cN'$, e.g., $\dbarbq'$ is the $\dbarb$-operator associated to $\cN'$. 

Let $X_{0}'$ be the section of $\cN'$ such that $\theta(X_{0}')=1$ and let $\varphi=\varphi_{X_{0},X_{0}'}$ be the vector bundle isomorphism of 
$T_{\C}M$ onto itself such that $\varphi$ is identity 
on $T_{1,0}\oplus T_{0,1}$ and $\varphi(X_{0})=X_{0}'$. By duality this defines a vector bundle isomorphism $\varphi^{t}$ from 
$\Lambda^{*}T^{*}_{\C}M$ onto itself. Then $\varphi^{t}$ induces an isomorphism from $\Lambda^{p,q}$ onto $\Lambda^{'p,q}$ and restricts to the 
identity on $\theta \wedge \Lambda^{*}T^{*}_{\C}M$. Thus, if $\eta \in C^{\infty}(M,\Lambda^{p,q})$ then $\varphi^{t}(\eta)$ is the component in 
$\Lambda^{'p,q}$ of $\eta$ with respect to the decomposition~(\ref{eq:CR-Lambda-pq-decomposition}) associated to $\cN'$. In fact, we can check that 
we have 
\begin{equation}
    \varphi^{t}(\eta)= \eta-\theta\wedge \iota_{X_{0}'}\eta \qquad \forall \eta \in C^{\infty}(M, \Lambda^{p,q}). 
     \label{eq:CR.varphit}
\end{equation}

\begin{lemma}\label{lem:CR.dependence-N}
    For $q=0,\ldots,n$ we have $\dbarbq'= \varphi^{t}\dbarbq (\varphi^{t})^{-1}$.
\end{lemma}
\begin{proof}
%
 Let $\eta \in C^{\infty}(M,\Lambda^{0,q})$ and let us compute $\dbarbq' [\varphi^{t}(\eta)]$. Thanks to~(\ref{eq:CR.varphit}) we have
 $d[\varphi^{t}(\eta)]=d\eta -d\theta \wedge 
 \iota_{X_{0}'}\eta +\theta\wedge d  \iota_{X_{0}'}\eta$. Moreover, we have 
 $\theta\wedge d  \iota_{X_{0}'}\eta=\theta \wedge \cL_{X_{0}'}\eta  - \theta \wedge \iota_{X_{0}'}d\eta$, so using~(\ref{eq:CR.dbarb}) we get 
 $\theta\wedge d  \iota_{X_{0}'}\eta=\theta \wedge \beta - \theta \wedge \iota_{X_{0}'}\dbarbq \eta -\theta\wedge \iota_{X_{0}'}{\partial}_{b;q}\eta$ 
 for some form $\beta$. Therefore, we see that  $d\varphi^{t}(\eta)$ is equal to
  \begin{multline}
      \theta \wedge \beta' + (\dbarbq \eta -\theta \wedge \iota_{X_{0}'}\dbarbq \eta) + 
     ({\partial}_{b;q} \eta -\theta \wedge \iota_{X_{0}'} {\partial}_{b;q} \eta) -d\theta \wedge  \iota_{X_{0}}\eta\\ 
     = \theta \wedge \beta' +  \varphi^{t}(\dbarbq \eta) + \varphi^{t}({\partial}_{b;q} \eta) - d\theta \wedge  \iota_{X_{0}}\eta,
      \label{eq:CR.functoriality-dbarbq}
 \end{multline}
 for some $\beta'$ in $C^{\infty}(M,\Lambda^{*}T^{*}_{\C}M)$. 
 
Since $T_{0,1}$ is integrable, for any sections $\overline{Z}$ and $\overline{W}$ of $T_{0,1}$ the Lie bracket $[\overline{Z},\overline{W}]$ is again a section of 
 $T_{0,1}$ and so we have 
 $d\theta(\overline{Z},\overline{W})=-\theta([\overline{Z},\overline{W}])=0$. This means that $d\theta$ has no component in $\Lambda^{`2,0}$. 
 Therefore, in~(\ref{eq:CR.functoriality-dbarbq}) the form $d\theta \wedge  \iota_{X_{0}}\eta$ cannot have a component in $\Lambda^{`0,q+1}$. 
 In view of the definition of $\dbarbq'$ it follows that $ \dbarbq' (\varphi^{t}(\eta))=  \varphi^{t}({\partial}_{b;q} \eta)$. Hence the lemma. 
%
%
 \end{proof}

%

Assume now that $M$ is endowed with a Hermitian metric $h$ on $T_{\C}M$ which commutes with complex conjugation and 
makes the splitting~(\ref{eq:CR-decomposition}) become orthogonal. 
Let $ \Box_{b;q}= \overline{\partial}^{*}_{b;q+1}\dbarbq +\overline{\partial}_{b;q-1}\overline{\partial}_{b;q}^{*}$ be the Kohn Laplacian and let 
$S_{b;q}$ be the Szeg\"o projection on $(0,q)$-forms, i.e., the orthogonal projection onto $\ker \Box_{b;q}$.
%

We also consider the orthogonal projections  $\Pi_{0}(\dbarbq)$ and $\Pi_{0}(\varthetabq)$ onto $\ker 
\dbarbq$ and $\ker \varthetabq=(\im \overline{\partial}_{b;q-1})^{\perp}$. In fact, as $\ker \dbarbq=\ker \Box_{b;q}\oplus \im 
\overline{\partial}_{b;q-1}$ we have $\Pi_{0}(\dbarbq)=S_{b;q}+1-\Pi_{0}(\varthetabq)$, that is, 
\begin{equation}
    S_{b;q}= \Pi_{0}(\dbarbq)+\Pi_{0}(\varthetabq)-1.
     \label{eq:Szego-projections-dbarb-varthetab}
\end{equation}
Let $N_{b;q}$ be the partial inverse of $\Box_{b;q}$, so that $N_{b;q}\Box_{b;q}=\Box_{b;q}N_{b;q}=1-S_{b;q}$. Then it can be shown
(see, e.g.,~\cite[pp.~170--172]{BG:CHM}) that we have
\begin{equation}
    \Pi_{0}(\dbarbq)=1-\overline{\partial}^{*}_{b;q+1}N_{b;q+1}\dbarbq, \qquad  
    \Pi_{0}(\varthetabq)=1-\overline{\partial}_{b;q-1}N_{b;q-1} \overline{\partial}_{b;q-1}^{*}.
     \label{eq:projections-dbarb-varthetab}
\end{equation}

The principal symbol of $\Box_{b;q}$ is invertible if, and only if, the condition $Y(q)$ holds at every point $x\in M$ (see~\cite[Sect.~21]{BG:CHM}, 
\cite[Sect.~3.5]{Po:MAMS1}). If we let $\kappa_{+}(x)$ and $\kappa_{-}(x)$ denote the number of positive and negative eigenvalues of $L_{\theta}$ at 
$x$, then the condition $Y(q)$ at  $x$ requires to have 
\begin{equation}
    q\not \in \{\kappa_{+}(x),\ldots,n-\kappa_{-}(x)\}\cup  \{\kappa_{-}(x),\ldots,n-\kappa_{+}(x)\}. 
\end{equation}

When the condition $Y(q)$ holds at every point the operator $\Box_{b;q}$  is hypoelliptic and 
admits a parametrix in $\pvdo^{-2}(M,\Lambda^{0,q})$ and then $S_{b;q}$ is a smoothing operator and $N_{b;q}$ is a \psivdo\ of 
order $-2$. Therefore, using~(\ref{eq:projections-dbarb-varthetab}) 
we see  that if the condition $Y(q+1)$ (resp.~$Y(q-1)$) holds 
everywhere then $\Pi_{0}(\dbarbq)$ (resp.~$\Pi_{0}(\varthetabq)$) is a \psivdo. 

Furthermore, in view of~(\ref{eq:Szego-projections-dbarb-varthetab}) we also see
 that if at every point the condition $Y(q)$ fails, but the conditions $Y(q-1)$ and 
$Y(q+1)$ hold, then the Szeg\"o projection $S_{b;q}$ is a zero'th order \psivdo\ projection. Notice that this may happen if, and only if, $M$ is 
$\kappa$-strictly pseudoconvex with $\kappa=q$ or $\kappa=n-q$.
 
Bearing all this in mind we have: 

\begin{theorem}\label{CR:Thm1}
    The following noncommutative residues are CR diffeomorphism invariants of $M$: \smallskip
    
    (i) $\Res \Pi_{0}(\dbarbq)$ when the condition $Y(q+1)$ holds everywhere;\smallskip 
    
    (ii) $\Res \Pi_{0}(\varthetabq)$ when the condition $Y(q-1)$ holds everywhere;\smallskip
    
    (iii) $\Res S_{b;\kappa}$ and $\Res S_{b;n-\kappa}$ when $M$ is $\kappa$-strictly pseudoconvex.\smallskip
    
\noindent In particular, they depend neither on the choice of the line bundle 
$\cN$, nor on that of the Hermitian metric $h$.
\end{theorem}
\begin{proof}
    Let us first show that the noncommutative residues in (i) and (ii) don't depend on the metric $h$. As the range of $\Pi_{0}(\dbarbq)$ 
    and the kernel of $\Pi_{0}(\dbarbq^{*})$ are $\ker \dbarbq$ and $(\ker \dbarbq^{*})^{\perp}=\im \dbarbq $ they don't depend on $h$. Therefore, if 
    the $Y(q+1)$ holds everywhere then $\Pi_{0}(\dbarbq)$ is a \psivdo\ projection whose range is independent of $h$, so the same its true for $\Res 
    \Pi_{0}(\dbarbq)$ by Lemma~\ref{prop:invariance-range-kernel}. Similarly, when the condition $Y(q-1)$ holds everywhere 
    the value of $\Res \Pi_{0}(\dbarbq^{*})$ is also independent of the choice of the Hermitian metric. 
    
    Next, let $\cN'$ be a supplement of $H$ in $TM$ and 
    let $h'$ be a Hermitian metric on $T_{\C}M$ which commutes with complex conjugation and 
    makes the splitting~(\ref{eq:CR-decomposition}) associated to $\cN'$ becomes orthogonal. We shall assign the superscript $'$ to objects associated 
    to the data $(\cN',h')$. 
    
    Let $X_{0}'$ be the section of $\cN'$ such that $\theta(X_{0}')=1$ and let $\varphi=\varphi_{X_{0},X_{0'}}$ be the vector bundle isomorphism of 
    $T_{\C}M$ onto itself such that $\varphi$ is identity on $T_{1,0}\oplus T_{0,1}$ and $\varphi(X_{0})=X_{0}'$. Since $\Pi_{0}(\dbarbq)$ 
    and $\Pi_{0}(\dbarbq^{*})$ don't depend on the choice of $h'$ we may assume that $h'=\varphi_{*}h$, so that $\varphi$ is a unitary isomorphism 
    from $(T_{\C}M,h)$ onto $(T_{\C},h')$ and $\varphi^{t}$ is a unitary vector bundle isomorphism from $\Lambda^{0,q}$ onto 
    $\Lambda^{'0,q}$. 
    
    Assume that the condition $Y(q+1)$ holds everywhere. Thanks to Lemma~\ref{lem:CR.dependence-N} we know that $\overline{\partial}_{b;q}'$ and 
    $\varphi^{t} \dbarbq (\varphi^{t})^{-1}$ agree. 
   Since $\varphi^{t}$ is unitary we also see that $\overline{\partial}_{b;q}^{*'}$ and $\varphi^{t} \dbarbq^{*} (\varphi^{t})^{-1}$ 
    too agree and so we have $\Box_{b;q}'=\varphi^{t} \dbarbq^{*} (\varphi^{t})^{-1}$. 
    Combining this with~(\ref{eq:projections-dbarb-varthetab}) we see that the projections
    $\Pi_{0}(\dbarbq')$ and $\varphi \Pi_{0}(\dbarbq) (\varphi^{t})^{-1}$ agree and so they have same noncommutative residue. 
    
    On the other hand, we have $\Res \varphi \Pi_{0}(\dbarbq) (\varphi^{t})^{-1}=\Res \Pi_{0}(\dbarbq)$, since $\Res$ is a trace. Hence $\Res 
    \Pi_{0}(\dbarbq')=\Res \Pi_{0}(\dbarbq)$, that is, the value of $\Res \Pi_{0}(\dbarbq)$ does not depend on $\cN$. In the 
    same way we can show that when the condition $Y(q-1)$ holds everywhere the residue $\Res \Pi_{0}(\dbarbq^{*})$ is independent of the choice made 
    for $\cN$. 

    Now, let $\phi:M\rightarrow M'$ be a CR diffeomorphism from $M$ onto a CR manifold $M'$. Let $\cN'$ be a supplement of $H$ in $TM$ and 
    let $h'$ be a Hermitian metric on $T_{\C}M'$ which commutes with complex conjugation and 
    makes the splitting~(\ref{eq:CR-decomposition}) of $T_{\C}M'$ 
    associated to $\cN'$ becomes orthogonal. We will assign the superscript $'$ to objects related to $M'$. 
    
    Since  the values of the noncommutative residues (i)--(ii) related to $M'$ are 
    independent of the data $(\cN, h)$, we may assume that $\cN=\phi_{*}\cN$ and $h'=\phi_{*}h$, so that $\phi$ gives rise to a unitary isomorphism from 
    $L^{2}(M,\Lambda^{0,q})$ onto $L^{2}(M',\Lambda^{'0,q})$.  As the fact that $\phi$ is a CR diffeomorphism implies that 
    $\phi_{*}\dbarbq=\dbarbq'$, we see that $\Pi_{0}(\dbarbq')=\phi_{*} \Pi_{0}(\dbarbq)$. Therefore, if the condition $Y(q+1)$ holds everywhere 
    then we have  $\Res\Pi_{0}(\dbarbq')= \Res \phi_{*} 
    \Pi_{0}(\dbarbq) = \Res \Pi_{0}(\dbarbq)$. 
    
    Similarly, when the condition $Y(q+1)$ holds everywhere we see that $\Res \Pi_{0}({\varthetabq}')$ and $\Res \Pi_{0}(\varthetabq)$ agree.  
Thus, the noncommutative residues (i) and (ii) are CR diffeomorphisms invariants.

    Finally, assume that $M$ is $\kappa$-strictly pseudoconvex and that $q=\kappa$ or $q=n-\kappa$. At every point of $M$ the condition $Y(q)$ fails, 
    but the condition $Y(q-1)$ and $Y(q+1)$ hold, so $S_{b;q}$ is a \psivdo\ projection and $\Res \Pi_{0}(\dbarbq)$ and $\Res \Pi_{0}(\varthetabq)$ 
    are CR diffeomorphism invariants.  
    Since~(\ref{eq:Szego-projections-dbarb-varthetab}) 
    implies that $\Res S_{b;q}=\Res \Pi_{0}(\dbarbq)+\Res 
    \Pi_{0}(\varthetabq)$, it follows that $\Res S_{b;q}$ too is an invariant of the CR diffeomorphism class of $M$. The proof is thus achieved.  
\end{proof}
%

Theorem~\ref{CR:Thm1} allows us to get CR invariants for CR manifolds that are not necessarily strictly pseudoconvex or have not necessarily a 
nondegenerate Levi form. However, specializing it to the strictly pseudoconvex case yields: 

\begin{theorem}
    Suppose that $M$ is a compact strictly pseudoconvex CR manifold. Then $\Res S_{b;k}$, $k=0,n$, and $\Res \Pi_{0}(\dbarbq)$, $q=1,\ldots,n-1$,  
are CR diffeomorphism invariants of $M$.  In particular, when $M$ is the boundary of a strictly pseudoconvex domain 
$D\subset \C^{n}$ they give rise to biholomorphism invariants of $D$.    
\end{theorem}

Finally, we can get further CR invariants by using the $\dbarb$-complex with coefficients in a CR holomorphic vector bundle as follows. 

A complex vector vector bundle $\cE$ over $M$ is a CR holomorphic vector bundle when there exists a patching of trivializations such that the transition maps are 
given by invertible 
matrices with CR function entries. For $q=0,\ldots, n$ let $\Lambda^{0,q}(\cE)=\Lambda^{0,q}\otimes \cE$. Then there exists a unique first order differential operator 
$\dbarbE: C^{\infty}(M,\Lambda^{0,*}(\cE))\rightarrow C^{\infty}(M, \Lambda^{0,*+1}(\cE))$ such that, for any local CR frame $e_{1},\ldots,e_{r}$ of $\cE$ and 
any local section $\omega=\sum \omega_{i}\otimes e_{i}$ of $\Lambda^{0,q}(\cE)$, we have 
\begin{equation}
    \dbarbE \omega = \sum_{i} (\dbarb \omega_{i})\otimes e_{i}. 
     \label{eq:CR.dbarbE}
\end{equation}
We have $\dbarbE^2=0$ and the Leibniz's rule holds, i.e., we have 
\begin{equation}
    \dbarbE(\eta \wedge \omega)= (\dbarb \eta)\wedge \omega +(-1)^q \eta\wedge \dbarbE \omega, 
\end{equation}
for any $(0,q)$-form $\eta$ and section $\omega$ of $\Lambda^{0,*}(\cE)$. Thus this yields a chain complex called the $\dbarb$-complex with coefficients in $\cE$. 

We equip $\cE$ with a Hermitian metric and let $\Box_{b,\cE}=\dbarbE^*\dbarbE + \dbarbE \dbarbE^*$ be the Kohn Laplacian with coefficients in $\cE$. It follows 
from~(\ref{eq:CR.dbarbE}) that in any CR trivialization $\Box_{b,\cE}$ has the same principal symbol as $\Box_{b}\otimes 1_{\cE}$, so its principal symbol is 
invertible if, and only if, the condition $Y(q)$ holds. Therefore, we can define the Szeg\"o projection $S_{b,\cE;q}$ and the projections $\Pi_{0}(\dbarbEq)$ and 
$\Pi_{0}(\dbarbEq^*)$ as before and we see that:\smallskip

- $\Pi_{0}(\dbarbEq)$ is a zero'th order projection under condition $Y(q+1)$;\smallskip

- $\Pi_{0}(\dbarbEq^*)$ is a zero'th order projection under condition $Y(q-1)$;\smallskip

- $S_{b,\cE;q}$  is a smoothing operator under the condition $Y(q)$, but it's a zero'th order projection when $M$ is $\kappa$-strictly pseudoconvex and $q=\kappa$ 
or $q=n-\kappa$.\smallskip

\noindent Then it is not difficult to modify the proof of Theorem~\ref{CR:Thm1} to get:

\begin{theorem}\label{CR:Thm2}
    The following noncommutative residues depend only on the CR diffeomorphism class of $M$ and of the CR holomorphic bundle isomorphism class of $\cE$: \smallskip
    
    (i) $\Res \Pi_{0}(\dbarbEq)$ when the condition $Y(q+1)$ holds everywhere;\smallskip 
    
    (ii) $\Res \Pi_{0}(\dbarbEq^*)$ when the condition $Y(q-1)$ holds everywhere;\smallskip
    
    (iii) $\Res S_{b,\cE;\kappa}$ and $\Res S_{b,\cE;n-\kappa}$ when $M$ is $\kappa$-strictly pseudoconvex.\smallskip
    
\noindent In particular, their values depend neither on the choice of the line bundle 
$\cN$, nor on that of the Hermitian metrics on $T_{\C}M$ and $\cE$.
\end{theorem}

\subsection{Invariance by deformation of the CR structure}
 We now look at the behavior of the CR invariants under deformations of the CR structure.  For sake of simplicity the results are proved 
 for the invariants of Theorem~\ref{CR:Thm1}, but they can be extended to the invariants of Theorem~\ref{CR:Thm2} with coefficients in a CR vector holomorphic 
 vector bundle $\cE$, provided that we consider deformations of the CR holomorphic structure of $\cE$ compatible with the deformation of the CR 
 structure of $M$.   
 
  First, let us look at what happens under deformations of the complex structure of $H$. Such a deformation is given by a smooth family 
  $(J_{t})_{t\in \R}\subset C^{\infty}(M,\End_{\R} H)$ such that for any $t \in \R$ we 
 have $J^{2}_{t}=-1$ and $T_{1,0,t}=\ker (J_{t}-i)$ is an integrable subbundle of $T_{\C}M$. 
 An important fact concerning such a deformations is: 
 
 \begin{lemma}\label{lem:CR.signature-Levi-form}
  Let $a \in M$. Then the signature of the Levi form $L_{\theta}$ in~(\ref{eq:CR.Levi-form}) at $a$
     is invariant under smooth deformations of the complex structure of $H$. Incidentally, for $q = 0,\ldots,n$ the condition $Y(q)$ at $a$ is 
     invariant under such deformations.
 \end{lemma}
 \begin{proof}
     First, let $J$ be the original complex structure of $H$. Then any section of $T_{1,0}=\ker (J-i)$ is of the form $X-iJX$ for some section $X$ of 
     $H$. Let $X$ and $Y$ be sections of $H$. Then the integrability condition on $T_{1,0}$ implies that we have 
     \begin{multline}
         0=\theta([X-iJX,Y-iJY])=-d\theta(X-iJX,Y-iJY)\\ = d\theta(JX,JY)-d\theta(X,Y)+i(d\theta(X,JY)+d\theta(JX,Y)),
     \end{multline}
 which gives $d\theta(JX,JY)=d\theta(X,Y)$ and $d\theta(JX,Y)=-d\theta(X,JY)$. Thus,  
  \begin{equation}
      L_{\theta}(X-iJX,Y-iY)=-id\theta(X-iJX,Y+iJY)=-2d\theta(JX,Y)-2id\theta(X,Y).
    \end{equation}
  Hence $X(a)-iJX(a)$ is in the kernel of $L_{\theta}(a)$ if, and only if, $X(a)$ and $JX(a)$ are in that of $d\theta(a)$. 
    
    Conversely, since $d\theta(a)(JX(a),Y(a))=-d\theta(X(a),JY(a))$ we see that $X(a)$ is in the kernel of $d\theta$ if, and only if, so is $JX(a)$. 
    It thus follows that the rank of 
    $d\theta(a)$ is twice that of $L_{\theta}(a)$. Incidentally, the latter is independent of the choice of $J$, hence is invariant under Heisenberg preserving 
    deformations of the CR structure. 
    
     Next, let $(J_{t})_{t \in \R}$ be a smooth family of complex structures on $H$ and for each $t$ let $L_{\theta,t}$ be the Levi form~(\ref{eq:CR.Levi-form}) 
     on $T_{1,0,t}=\ker (J_{t}-i)$. In order to show that the signature is independent of $t$ we only have to show that the number $\kappa(t)$ of its 
     negative eigenvalues is constant. Since the latter takes on integer values, it is actually enough to prove that it depends continuously on $t$. 
     
     Let $t_{0}\in \R$ and let $X_{1},\ldots,X_{2n}$ be a local frame of $H$ near $a$ such that $X_{n+j}=J_{t_{0}}X_{j}$ for $j=1,\ldots,n$. For $t\in 
     \R$ and $j=1,\ldots,n$ let $Z_{j,t}=X_{j}-iJ_{t}X_{j}$. Then $Z_{1,t}(a),\ldots,Z_{n,t}(a)$ depend 
     smoothly on $t$ and there exists $\delta>0$ such that $Z_{1,t}(a),\ldots,Z_{n,t}(a)$  form a basis of $T_{1,0;t;a}=\ker (J_{t}(a)-i)$ for 
     $|t_{0}-t|<\delta$.  For $t\in (t_{0}-\delta,t_{0}+\delta)$ let $A(t)=(-id\theta(a)(Z_{j,t}(a),\overline{Z_{k,t}}))_{1\leq j,k\leq 
     n}$ be the matrix of $L_{\theta,t}(a)$ with respect to this basis. This defines a smooth family of $n\times n$ Hermitian matrices. 
     
     Now, the rank of $L_{\theta,t}(a)$ and $\kappa(t)$ are respectively equal to the rank of $A(t)$ and to the number of its negative eigenvalues. 
     In particular, the rank of $A(t)$ is independent of $t$, say is equal to $r$. Then all the negative eigenvalues 
     $A(t)$ are contained in the 
     interval $[-\|A(t)\|,-\mu_{n-r}(A(t))]$, where $\mu_{n-r}(t)$ is the $(n-r)$'th eigenvalue of  $|A(t)|$, i.e., the absolute value of the 
     first non-zero eigenvalue of $A(t)$. Clearly $\|A(t)\|$ depends continuously on $t$, but the same is also true for $\mu_{n-r}(A(t))$
     as a consequence of the min-max principle. Therefore, for any $\delta'<\delta$ there exists real numbers $c_{1}$ and $c_{2}$ with $c_{2}<c_{1}<0$ 
     such that for $|t-t_{0}|\leq \delta'$ all the negative eigenvalues of $A(t)$ are contained in the interval $[c_{2},c_{1}]$. 
     
    For $t\in [t_{0}-\delta',t_{0}+\delta']$ let $E_{-}(A(t))$ be the negative eigenspace of $A(t)$ and let 
    $\Pi_{-}(A(t))$ be the orthogonal projection onto $E_{-}(A(t))$. Then $\kappa(t)$ is equal to $\dim E_{-}(A(t))=\Tr \Pi_{-}(A(t))$. Moreover, we 
    have the formula, 
    \begin{equation}
        \Pi_{-}(A(t))=\frac{1}{2i\pi}\int_{\Gamma}(\lambda-A(t))^{-1}d\lambda,
    \end{equation}
    where $\Gamma$ is any circle contained in the halfspace $\Re \lambda<0$ that  bounds an open disk containing $c_{2}$ and $c_{1}$. Since 
    $\Gamma$ is chosen independently of $t$ it follows from this that $ \Pi_{-}(A(t))$ is a continuous function of $t$ on 
    $[t_{0}-\delta',t_{0}+\delta']$. Hence $\kappa(t)=\Tr \Pi_{-}(A(t))$ is a continuous  function of $t$ near $t_{0}$. 
    
    All this shows that $\kappa(t)$ is a continuous function of $t$ on $\R$. As alluded to above this implies that the signature of $L_{\theta,t}$ is 
    independent of $t$. Hence the lemma. 
 \end{proof}
 
 \begin{proposition}\label{prop:CR.deformation} The invariants (i)--(iii) of Theorem~\ref{CR:Thm1} 
     are invariant under deformation of the complex structure of $H$.
 \end{proposition}
  \begin{proof}
      We will prove the result for $\Res \Pi_{0}(\dbarbq)$ only since the proofs for the other residues follow along similar lines.   
      
      Let $(J_{t})_{t \in \R}\subset C^{\infty}(M,\End_{\R}H)$ be a smooth family  of complex structures on $H$. 
      We can construct a smooth family of admissible Hermitian metric $h_{t}$ on $T_{\C}M$ as follows. 
 Let $g$ be a Riemannian metric on $H$ and let us extend it into a Hermitian metric on $H\otimes \C$ such that 
 \begin{equation}
     h(X_{1}+iY_{1},X_{2}+iX_{2})=g(X_{1},X_{2})+g(Y_{1},Y_{2})+i(g(Y_{1},X_{2})-g(X_{1},Y_{2})),
 \end{equation}
 for sections $X_{1},X_{2},Y_{1},Y_{2}$ of $H$. Notice that $h$ commutes with complex conjugation. 
 
 Let $X_{0}$ be the global section of $\cN$ such that 
 $\theta(X_{0})=1$. Then for any $t\in \R$ we get a Hermitian metric 
 $h_{t}$ on $T_{\C}M$ such that 
 \begin{equation}
     h_{t}(Z+\lambda X_{0},W+\mu X_{0})=g(Z,W)+g(J_{t}Z,J_{t}W)+\theta(X_{0})\theta(X_{0}')+\lambda \overline{\mu},
 \end{equation}
 for sections $Z,W$ of $H\otimes \C$ and functions $\lambda, \mu$ on $M$. This Hermitian  metric commutes with complex conjugation. Moreover, as $J_{t}$ is 
 unitary with respect to $h_{t|_{H}}$, the subbundles $T_{1,0;t}=\ker (J_{t}-i)$ and $T_{0,1;t}=\ker (J_{t}+i)$ are perpendicular with respect to 
 $h_{t}$, and so the splitting $T_{1,0;t}\oplus T_{0,1;t}\oplus (\cN\otimes \C)$ is orthogonal with respect to $h_{t}$.
 Therefore $(h_{t})_{t \in \R}$ is a smooth family of admissible Hermitian metrics on $T_{\C}M$. 
 
 We will use the subscript $t$ to denote operators related to the Hermitian metric $h_{t}$ and the CR structure defined 
 by $J_{t}$. In addition, we extend $J_{t}$ into a section of $\End T_{\C}M$ such that $J_{t}X_{0}=0$. Then $J_{t}$ commutes with its adjoint 
 $J_{t}^{*}=-J_{t}$ with respect to $h_{t}$, so the orthogonal projection $\pi_{0,1;t}$ onto $T_{0,1;t}=\ker (J_{t}+i)$ is  
 \begin{equation}
     \pi_{0,1,t}=\frac{1}{2i\pi}\int_{|\lambda+i|=1/2}(\lambda-J_{t})^{-1}d\lambda.
 \end{equation}
 In particular, $(\pi_{0,1;t})_{t\in \R}$ is a smooth family with values in $C^{\infty}(M,\End T_{\C}M)$. Henceforth, the family $(\pi_{t}^{0,q})\subset 
 C^{\infty}(M,\End T_{\C}^{*}M)$ of orthogonal projections onto $\Lambda^{0,q}_{t}$ depends smoothly on $t$. Therefore, 
 the operator $\dbarbqt=\pi^{0,q}_{t}\circ d$, its adjoint $\varthetabqt$ and the Kohn Laplacian $\Boxbqt$ give rise to smooth families of differential 
 operators, hence their principal symbols depend smoothly on~$t$. 
 
 Assume now that the condition $Y(q+1)$ holds everywhere. We know from Lemma~\ref{lem:CR.signature-Levi-form} 
 that this condition holds independently of $t$, so 
 the family of principal symbols $(\sigma_{2}(\Box_{b,t;q+1}))_{t\in \R}$ is a smooth family of invertible symbols. Therefore, by the results 
 of~\cite[Chap.~3]{Po:MAMS1} the family $(\sigma_{2}(\Box_{b,t;q+1})^{*-1})_{t \in \R}$ of the inverses is also a smooth family of symbols. 
 Since for any $t\in \R$ the principal symbol 
 of $N_{b,t;q+1}$ is $\sigma_{2}(\Box_{b,t;q+1})^{*-1}$, using~(\ref{eq:projections-dbarb-varthetab}) 
 we see that the principal symbol of $\Pi_{0}(\dbarbqt)$ depends smoothly on $t$. It then follows from 
 Lemma~\ref{lem:NCRPP.homotopy-invariance} 
 that $\Res \Pi_{0}(\dbarbqt)$ is independent of $t$. Hence $\Res \Pi_{0}(\dbarbq)$ is invariant under deformations of the complex structure of $H$. 
 \end{proof}

Finally, when $M$ is $\kappa$-strictly pseudoconvex we can deal with general deformations of the CR structure. 

\begin{proposition}\label{prop:CR.deformation-spc}
    Assume that $M$ is $\kappa$-strictly pseudoconvex. Then the invariants (i)--(iii) of Theorem~\ref{CR:Thm1} are invariant under deformations of the CR structure.
\end{proposition}
\begin{proof}
We shall prove the result for the invariant (i) only, since the proof is the same for the other invariants. 
    
Since $M$ is $\kappa$-strictly pseudoconvex a smooth deformation of the CR structure parametrized by a connected manifold $B$ consists of two 
pieces:\smallskip

- A smooth family of contact forms $(\theta_{\alpha})_{\alpha \in B}$, so that $d\theta_{\alpha}$ is nondegenerate on $H_{\alpha}=\ker 
\theta_{\alpha}$;\smallskip

- A family $(J_{\alpha})_{\alpha \in B}$  such that $J_{\alpha}$ is an (integrable) complex structures on $H_{\alpha}$ depending smoothly on 
$\alpha$.\smallskip 

Let $\alpha _{0}\in B$. As $(\theta_{\alpha})_{\alpha \in B}$ is a deformation of the contact structure of $M$, by a result of Gray~\cite[Sect.~5.1]{Gr:SGCS} 
there exists an open neighborhood $B'$ of $\alpha_{0}$ in $B$ and a smooth family $(\phi_{\alpha})_{\alpha \in B'}$ of diffeomorphisms of $M$ onto 
itself such that $\phi_{\alpha}^{*}H_{\alpha}=H_{\alpha_{0}}$. In addition, for $\alpha \in B$ and $\alpha \in B'$ let $T_{1,0;\alpha}=\ker (J_{\alpha}+i)$ and 
$T_{1,0;\alpha}'=\ker (\phi_{\alpha}^{*}J_{\alpha}+i)$. For $q \neq \kappa, n-\kappa$ we shall denote 
$\Res \Pi_{0;\alpha}(\overline{\partial}_{b,\alpha;p,q})$ (resp.~$\Res \Pi_{0;\alpha}'(\overline{\partial}_{b,\alpha;p,q}')$) the invariant (i) from 
Theorem~\ref{CR:Thm1} in bidegree $(0,q)$ associated to the CR structure defined by $T_{1,0;\alpha}$  (resp.~$T_{1,0;\alpha}'$). 

As $\phi_{\alpha}$ is a CR diffeomorphism from $(M, T_{1,0;\alpha}')$ onto $(M,T_{1,0;\alpha})$, by Theorem~\ref{CR:Thm1} we have 
$\Res \Pi_{0}(\dbarbq)_{\alpha}'=\Res \Pi_{0;\alpha}(\overline{\partial}_{b,\alpha;q})$. 
Observe also that $(\phi_{\alpha}^{*}J_{\alpha})$ is a smooth deformation of the 
complex structure of $H_{\alpha_{0}}$, so by Proposition~\ref{prop:CR.deformation} 
we have $\Res \Pi_{0;\alpha_{0}}(\overline{\partial}_{b,\alpha_{0};q})=\Res \Pi_{0}(\dbarbq)_{\alpha}'=\Res \Pi_{0;\alpha}(\overline{\partial}_{b,\alpha;q})$.  
Thus $\Res \Pi_{0;\alpha}(\overline{\partial}_{b,\alpha;q})$ is locally constant on $B$, hence constant since $B$  is connected. 
This shows that $\Res \Pi_{0}(\dbarbq)$ is invariant under deformations of the CR structure.
\end{proof}

\subsection{Computation of the invariants}\label{subsec:CR.computation}
Let us now make some comments about the computation the 
densities $c_{\Pi}(x)$ whose integrals yield the invariants $\Res \Pi$ from Theorems~\ref{CR:Thm1} and~\ref{CR:Thm1} (here $\Pi$ denotes any of 
the \psivdo\ projection involved in these Theorems). 

As explained in Introduction the computation of the densities $c_{\Pi}(x)$ is interesting even if $\Res \Pi$ may vanish, because it could provide us 
with geometric information about the logarithmic singularity of the kernel of the geometric projection $\Pi$. However, the direct computation of 
$c_{\Pi}(x)$ in local coordinates is rather involved: it amounts to determine the symbol of degree 
$-(2n+2)$ of $\Pi$, so that we have more and more terms to compute as the dimension increases. Therefore, we need additional tools to deal with the computation. 

%

When the bundle $\cE$ is trivial and the CR manifold $M$ is strictly pseudoconvex and endowed with the Levi metric defined by a pseudohermitian 
contact form  $\theta$,
we can extend the arguments of~\cite{BGS:HECRM} to show that 
the densities $c_{\Pi}(x)$ are of the form $\tilde{c}_{\Pi}(x)d\theta^{n}\wedge \theta$, where $\tilde{c}_{\Pi}(x)$ is a local pseudohermitian 
invariant of weight $n+1$. This means that  $\tilde{c}_{\Pi}(x)$ is a universal polynomial in complete contractions of the covariant derivatives of the curvature and 
torsion tensors of the 
Tanaka-Webster connection and the polynomial is homogeneous of degree $-(n+1)$ under scalings $\theta\rightarrow \lambda\theta$, $\lambda>0$, of the 
pseudohermitian contact form. Thus the residues $\Res \Pi$ are \emph{geometric} global CR invariants. In dimension 3 there are no non-zero such 
invariants (see~\cite{BHR:DLEICRMD3}), but to date there no known obstruction to the existence of global geometric CR invariant.

In conformal geometry the conjecture of Deson-Swimmer~\cite{DS:GCCAAD}, partially proved by Alexakis~(\cite{Al:DGCI1}, \cite{Al:DGCI2}), 
predicts the local form of the Riemannian 
invariants whose integrals yield global conformal invariants. It would be very interesting to prove an analogue of this conjecture in CR geometry, 
but to date it is not even clear what could be the conjecture, so we cannot use it to predict the form of the densities $c_{\Pi}(x)$. However, the 
computation of some $c_{\Pi}(x)$ by other means would certainly shed some light on some of the pseudohermtian invariants that should enter in the 
conjecture.  

On the other hand, in the case of the Szeg\"o projection $S_{b,0}$ on functions the density $c_{S_{b,0}}(x)$ is not a CR invariant, 
but it transforms conformally under the conformal 
changes $\theta \rightarrow e^{2f}\theta$ of pseudohermitian contact forms that come from CR pluriharmonic functions $f$, i.e., functions that are locally  
real parts of CR functions. Therefore, it would be natural to try to extend the CR invariant theory of~\cite{BEG:ITCCRG} and~\cite{Hi:CBILSBK} to 
deal with this class of invariants and to get information about the logarithmic singularity of the Szeg\"o kernel. It seems that Hirachi~\cite{Hi:PC} has made 
recent progress in this direction.

We could like to suggest another approach which comes in naturally with the framework of the paper and would allows us to deal with invariants with 
coefficients in CR holomorphic vector bundles as well. Namely,  it would be natural to make use of a version of rescaling of 
Getzler~\cite{Ge:SPLASIT} to simplify the computation of the noncommutative residues densities. The latter is a powerful trick which, 
by taking into account the supersymmetry of the Dirac operator, allows us to get a short proof of
the small time convergence to the Atiyah-Singer integrand of the local supertrace of 
the heat kernel 
of the square of the Dirac operator. This bypasses the invariant theory of~\cite{ABP:OHEIT} and \cite{Gi:ITHEASIT}
and provides us with a purely analytical proof 
of the Atiyah-Singer index theorem for Dirac operators. 

Let $M$ be a strictly pseudoconvex CR manifold endowed with the Levi metric defined by pseudohermitian contact form $\theta$ and let $\cE$ be 
a Hermitian CR vector bundle over $M$.  It is believed that implementing a version of Getzler's rescaling into the Heisenberg calculus would allow us to compute the 
density 
\begin{equation}
    \Str_{\Lambda^{0,*}(\cE)} c_{S_{b,\cE}}(x)=c_{S_{b,\cE;0}}(x)+(-1)^{n}\Tr_{\Lambda^{0,n}(\cE)}c_{S_{b,\cE;n}}(x),
\end{equation}
where $\Str:=\sum_{q=0}^{n}(-1)^{q}\Tr_{\Lambda^{0,q}(\cE)}$ is the supertrace on $\Lambda^{0,*}(\cE)$ and $S_{b,\cE}$ denotes the Szeg\"o projection acting 
on all sections of $\Lambda^{0,*}(\cE)$. 

It should be apparent from~\cite{Ge:ADECRM} and~\cite{Po:CMP1} that Getzler's rescaling techniques could be used in the setting of the 
Heisenberg calculus. The upshot is that the Getzler's rescaling would yield near any point of the manifold a refinement of the filtration of 
the Heisenberg calculus, so that determining $ \Str c_{S_{b,\cE}}(x)$ would boil down to computing  the second subleading symbol of $S_{b,\cE}$ with 
respect to this new filtration. This would be 
\emph{infinitely} better than to have to compute the symbol of order $-(2n+2)$ in the usual sense of the Heisenberg calculus. 
We expect to carry out the explicit calculation in a future paper. 
%
%

\section{Invariants of generalized Szeg\"o projections}\label{sec.Toeplitz}
Let $(M^{2n+1},H)$ be an orientable contact manifold, i.e., a Heisenberg manifold admitting a real 1-form $\theta$, called contact form, 
such that $\theta$ annihilates $H$ and $d\theta_{|_{H}}$ is nondegenerate. 
Given a contact form $\theta$ on $M$ we let $X_{0}$ be the Reeb 
vector field of $\theta$, i.e., the unique vector field $X_{0}$ such that $\iota_{X_{0}}\theta=1$ and $\iota_{X_{0}}d\theta=0$. 

In addition, we let $J$ be an 
 almost complex structure on $H$  which is \emph{calibrated} in the sense that $J$ preserves $d\theta_{|H}$ and we have $d\theta(X,JX)>0$ for any 
 non-vanishing section $X$ of $H$. Extending $J$ to $TM$ by requiring to have  $JX_{0}=0$, we then can equip $TM$ with the Riemannian metric 
 $g_{\theta, J}:=d\theta(.,J.)+\theta^{2}$.

In this context Szeg\"o projections have been defined by Boutet de Monvel and Guillemin in~\cite{BG:STTO} as FIO's with complex phase. 
This construction has been further generalized by Epstein-Melrose~\cite{EM:HAITH}  as follows. 

Let $\bH^{2n+1}$ be the Heisenberg group of dimension $2n+1$ consisting of $\R^{2n+1}$ together with the 
group law,
\begin{equation}
     x.y=(x_{0}+y_{0}+\frac{1}{2}\sum_{1\leq j\leq n}(x_{n+j}y_{j}-x_{j}y_{n+j}),x_{1}+y_{1},\ldots,x_{2n}+y_{2n}).
\end{equation}

Let $\theta^{0}=dx_{0}+\frac{1}{2}\sum_{j=1}^{n}(x_{j}dx_{n+j}-x_{n+j}dx_{j})$ be the standard left-invariant contact form of $\bH^{2n+1}$. Its Reeb 
vector field is $X^{0}_{0}=\frac{\partial}{\partial x_{0}}$  and if 
for $j=1,\ldots,n$  we let $X^{0}_{j}=\frac{\partial}{\partial x_{j}}+\frac{1}{2}x_{n+j}\frac{\partial}{\partial 
    x_{0}}$ and $X^{0}_{n+j}=\frac{\partial}{\partial x_{n+j}}-\frac{1}{2}x_{j}\frac{\partial}{\partial 
    x_{0}}$ then $X_{1}^{0},\ldots,X_{2n}^{0}$ form a left-invariant frame of $H^{0}=\ker \theta^{0}$. Note that for $j,k=1,\ldots,n$ and $k\neq j$ we have 
 \begin{equation}
     [X^{0}_{j},X^{0}_{n+k}]=-\delta_{jk}X^{0}_{0}, \qquad [X^{0}_{0},X^{0}_{j}]=[X^{0}_{j},X^{0}_{k}]=[X^{0}_{n+j},X^{0}_{n+k}]=0.
     \label{eq:IGSP.Heisenberg-relations}
 \end{equation}
%
    
The standard CR structure of $\bH^{2n+1}$ is given by the complex structure $J^{0}$ on $H^{0}$ such that $J^{0}X_{j}^{0}=X_{n+j}^{0}$ and 
 $J^{0}X_{n+j}=-X_{j}$.  It follows from~(\ref{eq:IGSP.Heisenberg-relations}) 
 that $J^{0}$ is calibrated  with respect to $\theta^{0}$ and that 
 $X_{0}^{0}, X_{1}^{0},\ldots,X_{2n}^{0}$ form an orthonormal frame of 
 $T\bH^{2n+1}$ with respect to the metric $g_{\theta^{0},J^{0}}$. 
 
 The scalar Kohn Laplacian  on $\bH^{2n+1}$ is equal to
\begin{equation}
     \Box_{b,0}^{0}=-\frac{1}{2}((X_{1}^{0})^{2}+\ldots+(X_{2n}^{0})^{2})+i\frac{n}{2}X_{0}^{0}.
\end{equation}
For $\lambda \in \C$ the operator $-\frac{1}{2}((X^{0}_{1})^{2}+\ldots+(X^{0}_{2n})^{2})+i \lambda X^{0}_{0}$ is invertible if, and only if,  we have $\lambda \not \in  
\pm(\frac{n}{2}+\N)$ (see~\cite{FS:EDdbarbCAHG}, \cite{BG:CHM}). For $k=0,1,\ldots$ the orthogonal projection 
$\Pi_{0}(\Box_{b}+ikX^{0}_{0})$ onto the kernel of $\Box_{b}+ikX^{0}_{0}$ is a left-invariant homogeneous  \psivdo\ of order~$0$ 
(see~\cite[Thm.~6.61]{BG:CHM}).  We then let $s_{k}^{0}\in S_{0}((\fh^{2n+1})^{*})$ 
denote its symbol, so that  we have $\Pi_{0}(\Box_{b}+ikX^{0}_{0})=s_{k}^{0}(-iX^{0})$. 
%
%
%
%
%
%

Now, since the existence of a contact structure implies that the Levi form~(\ref{eq:NCRP.Levi-form}) of $(M,H)$ is everywhere nondegenerate,
the tangent Lie group bundle $GM$ is a fiber bundle with typical fiber 
$\bH^{2n+1}$ (see~\cite{Po:Pacific1}). A local trivialization near a given point $a \in M$ is obtained as follows. 

Let $X_{1},\ldots,X_{2n}$ be a local orthonormal frame of $H$ on an open neighborhood $U$ of $a$ and which is admissible in the sense 
that $X_{n+j}=JX_{j}$ for $j=1,\ldots,n$. In addition, let $\underline{X_{0}(a)}$ denote the 
class of $X_{0}(a)$ in $T_{a}M/H_{a}$. Then as shown in~\cite{Po:Pacific1} 
the map $\phi_{X,a}:(T_{a}M/H_{a})\oplus H_{a}\rightarrow \R^{2n+1}$ such that 
\begin{equation}
    \phi_{X,a}(x_{0}\underline{X_{0}(a)}+x_{1}X_{1}(a)+\ldots + 
    x_{2n}X_{2n}(a))=(x_{0},\ldots,x_{2n}), \qquad x_{j}\in \R, 
\end{equation}
gives rise to a Lie group isomorphism from $G_{a}M$ onto $\bH^{2n+1}$. In fact, as $\phi_{X,a}$ depends smoothly on $a$ we get a  
fiber bundle trivialization of $GM|_{U}\simeq U\times \bH^{2n+1}$. 

For $j=0,\ldots,2n$ let $X_{j}^{a}$ be the model vector field of $X_{j}$ at $a$ as defined in~\cite{Po:Pacific1}. This is 
the unique left-invariant vector field on $G_{a}M$ which, in the coordinates provided by $\phi_{X,a}$, agrees with $\frac{\partial}{\partial x_{j}}$ at 
$x=0$. Therefore, we have $X_{j}^{a}=\phi_{X,a}^{*}X_{j}^{0}$ and so we get
\(
    \phi_{X,a}^{*}\Box_{b}^{0}=-\frac{1}{2}((X_{1}^{a})^{2}+\ldots+(X_{2n}^{a})^{2})+i\frac{n}{2}X_{0}^{a}. 
\)

If $\tilde{X}_{1},\ldots,\tilde{X}_{2n}$ is another admissible orthonormal frame of $H$ near $a$, then we pass from 
$(\tilde{X}_{1}^{a},\ldots,\tilde{X}_{2n}^{a})$ to $(X_{1}^{a},\ldots,X_{2n}^{a})$ by an orthogonal linear transformation, which leaves the 
expression $(X_{1}^{a})^{2}+\ldots+(X_{2n}^{a})^{2}$ unchanged. Therefore, the differential operator $\Box_{b}^{a}:=\phi_{X,a}^{*}\Box_{b}^{0}$ makes sense 
independently of the choice of the admissible frame $X_{1},\ldots,X_{2n}$ near $a$. 

On the other hand, as $\phi_{X,a}$ induces a unitary transformation from $L^{2}(G_{a}M)$ onto $L^{2}(\bH^{2n+1})$ we have 
$\Pi_{0}(\Box_{b}^{a}+ikX_{0}^{a})=\Pi_{0}(\phi_{X,a}^{*}(\Box_{b}^{0}+ikX_{0}^{0}))=\phi_{X,a}^{*}\Pi_{0}(\Box_{b}^{0}+ikX_{0}^{0}))$. 
 Hence $\Pi_{0}(\Box_{b}^{a}+ikX_{0}^{a})$ is a zero'th order left-invariant homogeneous \psivdo\ on $G_{a}M$ with symbol 
$s_{k}^{a}(\xi)=\phi_{X,a}^{*}s_{k}^{0}(\xi)=s_{k}^{0}((\phi_{X,a}^{-1})^{t}\xi)$. 
 In fact, since $\phi_{X,a}$ depends smoothly on $a$ we obtain:

\begin{proposition}[{\cite[Chap.~6]{EM:HAITH}}]
    For $k=0,1,\ldots$ there is a uniquely defined symbol $s_{k}\in S_{0}(\fg^{*}M)$ such that, for any admissible orthonormal frame $X_{1},\ldots,X_{d}$ of $H$ 
    near a 
    point $a \in M$, we have $s_{k}(a,\xi)=\phi_{X,a}^{*}s_{k}^{0}(\xi)$ for any $(a,\xi)\in \fg^{*}M\setminus 0$.
\end{proposition}

We call $s_{k}$ the \emph{Szeg\"o symbol at level} $k$. This definition \emph{a priori} depends on the contact form $\theta$ and the almost complex 
structure $J$. As we shall now see changing $\theta$ or $J$ has minor effects on $s_{k}$, but first we need the following. 

\begin{lemma}\label{lem:calibrated-almost-complex-structures}
    The space of calibrated almost complex structures on $H$ is path-connected.
\end{lemma}
\begin{proof}
    Let $J'\in C^{\infty}(M,\End H)$ be a calibrated almost complex structure on $H$ and set $g_{0}=g_{\theta,J}$ and 
    $g_{1}=g_{\theta,J'}$.  In the sequel the transpose superscript refers to tranposition with respect to $g_{0}$. 
    For any sections $X$ and $Y$ of $H$ we have $g_{1}(X,Y)=d\theta(X,J'Y)=g_{0}(JX,J'Y)=g_{0}(J^{'t}JX,Y)$. In particular, we see that 
    $A_{1}:=J^{'t}J$ is a  symmetric and positive definite section of $\End H$. Furthermore, as $J'$ preserves $d\theta_{|H}$ we have 
    $g_{1}(X,Y)=d\theta(X,J'Y)=-d\theta(J'X,Y)=-g_{0}(JJ'X,Y)$. Hence $J^{'t}J'=-JJ'$, which gives $J^{'t}=JJ'J$.
    
    For $t\in [0,1]$ let $g_{t}=(1-t)g_{0}+tg_{1}$. This is a smooth family of Riemannian metrics on $TM$. On $H$ we have $g_{t}(X,Y)=g_{0}(A_{t}X,Y)$ 
    with $A_{t}=(1-t).1+tJ^{'t}J$ and we can write $d\theta(X,Y)=g_{0}(JX,Y)=g_{t}(B_{t}X,Y)$ with $B_{t}=A_{t}^{-1}J$. Notice that $B_{t}$ is 
    antisymmetric with respect to $g_{t}$. Moreover, for $t=0$ we have 
    $B_{0}=J$ and since $J^{'t}=JJ'J$ for $t=1$ we have $B_{1}=(J^{'t}J)^{-1}J=JJ^{'t}J=J'$. 
    
    Since $B_{t}$ is antisymmetric with respect to $g_{t}$ its modulus and its phase with respect to $g_{t}$ are 
    $|B_{t}|_{t}=\sqrt{-B_{t}^{2}}$ and $J_{t}=B_{t}(\sqrt{-B_{t}^{2}})^{-1}$. Thus $(J_{t})_{t\in [0,1]}$ is a smooth path in $C^{\infty}(M,\End H)$ 
    such that $J_{t}$ is orthogonal with respect to $g_{t}$ for any $t\in [0,1]$. Notice that $B_{0}=J$ is already orthogonal with respect  to 
    $g_{0}$, so we have $J_{0}=J$. Similarly, we have $J_{1}=J'$. As $B_{t}$ is antisymmetric with respect to 
    $g_{t}$ the same is true for $J_{t}$. Together with the orthogonality this implies that we have $J_{t}^{-1}=-J_{t}$, i.e., $J_{t}$ is an almost complex 
    structure on $H$. Moreover, for sections $X$ and $Y$ of $H$ with $X$ non-vanishing, we have 
    $d\theta(X,J_{t}X)=g_{t}(B_{t}X,J_{t}X)=g_{t}(J_{t}^{-1}B_{t}X,X)=g_{t}(|B_{t}|_{t}X,X)>0$ and $d\theta(J_{t}X,J_{t}Y)=g_{t}(B_{t}J_{t}X,J_{t}Y)=
    g_{t}(J_{t}B_{t}X,J_{t}Y)=g_{t}(B_{t}X,Y)=d\theta(X,Y)$. Therefore $(J_{t})_{t\in [0,1]}$ is a smooth path of calibrated almost complex 
    structures on $H$ connecting $J$ to $J'$. Hence the lemma.
\end{proof}

Granted this we can now prove:

\begin{lemma}[{\cite[Chap.~6]{EM:HAITH}}]\label{lem:IGSP.properties-sk}
    (i) $s_{k}$ is invariant under conformal changes of contact form.\smallskip
    
    (ii) The change $(\theta, J)\rightarrow (-\theta, -J)$ transforms $s_{k}$ into $s_{k}(x,-\xi)$.\smallskip
    
    (iii) $s_{k}$ depends on $J$ only up to homotopy of idempotents in $S_{0}(\fg^{*}M)$.
 \end{lemma}
\begin{proof}
  Throughout the proof we let  $X_{1},\ldots,X_{2n}$ be an admissible orthonormal frame of $H$ near a point $a\in M$.
 
  Let $\theta'$ be a contact form which is conformal to $\theta$, that is,  $\theta'=e^{-2f}\theta$ with $f \in C^{\infty}(M,\R)$, 
  and let $s_{k}'$ be the Szeg\"o symbol at level $k$ 
   with respect to $\theta'$ and $J$. For $j=1,\ldots,2n$ let $X_{j}'=e^{f}X_{j}$. Then $X_{1}',\ldots,X_{2n}'$ is an admissible orthonormal frame of $H$ 
   with respect to $g_{\theta', J|_{H}}=e^{2f}d\theta(.,J.)$. Moreover, as the Reeb 
   vector field of $\theta'$ is such that $X_{0}'=e^{2f}X_{0} \bmod H$, we have $\underline{X_{0}'}(a)=e^{2f(a)}\underline{X_{0}}(a)$. 
   Thus, 
   \begin{multline}
       \phi_{X',a}^{-1}(x_{0},\ldots,x_{2n})=x_{0}\underline{X_{0}'}(a)+x_{1}X_{1}'(a)+\ldots+x_{2n}X_{2n}'(a)\\ = 
       x_{0}\lambda^{2}\underline{X_{0}}(a)+x_{1}\lambda X_{1}(a)+\ldots+x_{2n}\lambda X_{2n}(a))  
       =\phi_{X,\alpha}^{-1}\circ \delta_{\lambda}(x_{0},\ldots,x_{2n})
        \label{eq:IGSP.indetificationGaM-H}
   \end{multline}
   where $\lambda=e^{f(a)}$ and $\delta_{\lambda}(x)=\lambda.x$ for any $x \in \bH^{2n+1}$. 
   
   On the other hand, as $s_{k}^{0}$ is homogeneous of degree $0$ we have $\delta_{\lambda *}s_{k}^{0}=s_{k}^{0}$. Therefore, we get
   $s_{k}'(a,.)=\phi_{X',a}^{*}s_{k}^{0}=\phi_{X,a}^{*}\delta_{\lambda *}s_{k}^{0}=\phi_{X,a}^{*}s_{k}^{0}=s_{k}(a,.)$.
   Hence $s_{k}$ is invariant under conformal changes of contact form.
   
   Let $s_{k}'$ be the Szeg\"o symbol at level $k$ with respect to $-\theta$ and $-J$. Define $X_{0}'=-X_{0}$ and for $j=1,\ldots,n$ 
   let $X'_{j}=X_{j}$ and $X_{n+j}'=-X_{n+j}$. Then $X_{0}'$ is the Reeb vector field of $-\theta$ and $X_{1}',\ldots,X_{2n}'$ form an admissible 
   orthonormal frame with respect to $g_{-\theta,-J}$. Moreover, we have $\phi_{X',a}=\tau\circ \phi_{X,a}$, where we have let
   $\tau(x)=(-x_{0},x_{1},\ldots,x_{n},-x_{n+1},\ldots,x_{2n})$. Hence $s_{k}'(a,.)=\phi_{X',a}s_{k}^{0}=\phi_{X,a}^{*}\tau^{*}s_{k}^{0}$. 
     
   We have $ \tau^{*}[s_{k}^{0}(-iX^{0})]=\tau^{*}\Pi_{0}(\Box_{b}^{0}+ikX_{0}^{0})= \Pi_{0}(\tau^{*}(\Box_{b}^{0}+ikX_{0}^{0}))$, since 
   the action of $\tau$ on $L^{2}(\bH^{2n+1})$ is unitary.  Moreover, as $\tau^{*}X_{j}^{0}$ is equal to $X_{j}^{0}$ if $j=1 ,\ldots,n$ and to
   $-X_{j}^{0}$ otherwise, we see that  
   $\tau^{*}(\Box_{b}^{0}+ikX_{0}^{0})$ is equal to 
\begin{equation}
       -\frac{1}{2}((X_{1}^{0})^{2}+\ldots+(X_{2n}^{0})^{2})-i(\frac{n}{2}+k)X_{0}^{0}= 
       (\Box_{b}^{0}+ikX_{0}^{0})^{t}.
\end{equation}
Thus $\tau^{*}[s_{k}^{0}(-iX^{0})]=\Pi_{0}((\Box_{b}^{0}+ikX_{0}^{0})^{t})=[\Pi_{0}(\Box_{b}^{0}+ikX_{0}^{0})]^{t}=s_{k}^{0}(-iX^{0})^{t}$, 
using~(\ref{eq:HC.symbols-transpose-adjoints}) we obtain $(\tau^{*}s_{k}^{0})(\xi)=s_{k}^{0}(-\xi)$. Hence we have
   $s_{k}'(a,\xi)=(\phi_{X,a}^{*}\tau^{*}s_{k}^{0})(\xi)=\phi_{X,a}^{*}(s_{k}^{0}(-\xi))=s_{k}(a,-\xi)$.

   
  Let $J'$ be another almost complex structure on $H$ calibrated with respect to $\theta$ and let $s_{k}'$ be the Szeg\"o symbol at level $k$ with 
    respect to $\theta$ and $J'$. Then by Lemma~\ref{lem:calibrated-almost-complex-structures} there exists a smooth path $(J_{t})_{0\leq t\leq 0}$ 
    of calibrated almost complex structures such that $J_{0}=J$ and $J_{1}=J'$. 
      For $t\in [0,1]$ let $X_{j,t}=X_{j}$ if $j=0,1,\ldots,n$ and $X_{j,t}=J_{t}X_{j}$ otherwise. Then $X_{1,t},\ldots,X_{2n,t}$ is an admissible 
   orthonormal frame of $H$ with respect to $g_{\theta,J_{t}}$ and the isomorphism $\phi_{X_{t},a}:G_{a}M\rightarrow \bH^{2n+1}$ 
   depends smoothly on $a$ and $t$. Therefore, $s_{k,t}(a,\xi)=\phi_{X_{t},a}^{*}s_{k}^{0}(\xi)$ a smooth path  in $S_{0}((\fh^{2n+1})^{*})$ 
         connecting $s_{k}$ to $s_{k}'$.  The proof is now complete.
\end{proof}
 
From now on we let $\cE$ be a Hermitian vector bundle over $M$. 
\begin{definition}[{\cite[Chap.~6]{EM:HAITH}}]\label{def:IGSP-definition}
    For $k=0,1,\ldots$ a generalized Szeg\"o projection at level $k$ is a \psivdo\ projection $S_{k}\in \pvdo^{0}(M,\cE)$ with principal symbol 
    $s_{k}\otimes \op{id}_{\cE}$.
\end{definition}

Generalized Szeg\"o projections at level $k$ always exist (see~\cite[Prop.~6]{EM:HAITH} and Lemma~\ref{lem:K-Theory.symbol-projection}). Moreover, when 
$k=0$  and $\cE$ is the trivial line bundle the above definition allows us to recover the Szeg\"o projections of~\cite{BG:STTO}, for we have: 

\begin{lemma}\label{lem:Contact.example}
 Let $S:C^{\infty}(M)\rightarrow C^{\infty}(M)$ be a 
Szeg\"o projection in the sense of~\cite{BG:STTO}. Then $S$ is a generalized Szeg\"o projection at level $0$.
\end{lemma}
\begin{proof}
    We saw in Section~\ref{sec.Heisenberg} that $S$ is a zero'th order \psivdo. Moreover, if $q(x,y)$ is the complex 
    phase of $S$ then it follows from~(\ref{eq:NRHC.Taylor-complex-phase}) 
    that at a point $a\in M$ the model operator $S^{a}$ of $S$ in the sense of~\cite[Def.~3.2.7]{Po:MAMS1} is a Szeg\"o projection, 
    whose complex phase is given by the 
    leading term at $x=a$ in~(\ref{eq:NRHC.Taylor-complex-phase}). In particular, under the identification $G_{a}M\simeq \bH^{2n+1}$ provided by a map 
    $\phi_{X,a}$ as in~(\ref{eq:IGSP.indetificationGaM-H}) we 
    see that $(\phi_{X,a})_{*}S^{a}$ is a Szeg\"o projection on $\bH^{2n+1}$. In fact, as $(\phi_{X,a})_{*}S^{a}$ is left-invariant and homogeneous this is the 
    Szeg\"o projection 
    $\Pi_{0}(\Box_{b})$ considered above, so that $(\phi_{X,a})_{*}S^{a}$ has symbol $s_{0}^{0}$. Since by definition $S^{a}$ has symbol $\sigma_{0}(S)(a,.)$ 
    we see that  $\sigma_{0}(S)=s_{0}$. Hence $S$ is a generalized Szeg\"o projection at level $0$. 
\end{proof}


In particular, when $M$ is strictly pseudoconvex the Szeg\"o projection $S_{b,0}$ is a generalized Szeg\"o projection at level $0$. 

If $S_{k}$ and $S_{k}'$  are generalized Szeg\"o projection at level $k$ in $\pvdo^{0}(M,\cE)$ then they have same principal, so by 
Lemma~\ref{lem:NCRPP.homotopy-invariance} they have same noncommutative residue. We then let 
\begin{equation}
    L_{k}(M,\cE):=\Res S_{k}.
\end{equation}

Recall that the $K$-group $K^{0}(M)$ can be described as the group of formal differences of stable homotopy classes of (smooth) vector bundles over $M$, 
where a stable homotopy between vector bundles $\cE_{1}$ and $\cE_{2}$ is given by an auxiliary vector bundle $\cF$ and a vector bundle isomorphism 
$\phi:\cE_{1}\oplus \cF \simeq \cE_{2}\oplus \cF$. Then we have:

\begin{theorem}\label{thm:IGSP.main}
$L_{k}(M,\cE)$ depends only on the Heisenberg diffeomorphism class of $M$ and on the $K$-theory class of $\cE$, hence is invariant under deformations of the 
contact structure. In particular, $L_{k}(M,\cE)$ depends neither on the contact form $\theta$, nor on the almost complex structure $J$.
\end{theorem}
\begin{proof}
    Throughout the proof we let $S_{k}\in \pvdo^{0}(M,\cE)$  be generalized Szeg\"o projection at level $k$, so that  $L_{k}(M,\cE)=\Res S_{k}$.
    
    Let us first show that $L_{k}(M,\cE)$  is independent from $\theta$ and $J$. To this end let $\theta'$ be 
    a contact form on $M$, let $J'$ be an almost complex structure on $H$ calibrated with respect to $\theta'$ and let $S_{k}'\in \pvdo^{0}(M,\cE)$ be a 
    generalized Szeg\"o projection at level $k$ with respect to $\theta'$ and $J'$. 
    
    If $\theta$ and $\theta'$ are in the same conformal class, then by Lemma~\ref{lem:IGSP.properties-sk}
    the principal symbols of $S_{k}$ and $S_{k}'$ are 
    homotopic, so by Lemma~\ref{lem:NCRPP.homotopy-invariance}  we have $\Res S_{k}'=\Res S_{k}$. 
    
    Let $\tau:\cE \rightarrow \cE^{*}$ be the antilinear isomorphism provided by the Hermitian metric of $\cE$ and define 
    $S_{k}''=\tau^{-1}S_{k}^{t}\tau$. Then  $S_{k}''$  is a \psivdo\ projection and by~(\ref{eq:HC.symbols-transpose-adjoints}) its principal symbol is
    $\tau^{-1}(s_{k}(x,-\xi)\otimes \op{id}_{\cE^{*}})\tau=s_{k}(x,-\xi)\otimes \op{id}_{\cE}$. Therefore, by Lemma~\ref{lem:IGSP.properties-sk}
    this is a generalized Szeg\"o projection at level $k$ with respect to $-\theta$ and 
    $-J$. Thus, if $\theta'$ is not in the conformal class of $\theta$ then it is in that of $-\theta$ and as above we get $\Res 
    S_{k}'=\Res S_{k}''$. As $\Res S_{k}''=\Res S_{k}^{t}$ and by Lemma~\ref{NCRP.real-valued} we have $\Res S_{k}^{t}=\Res S_{k}$, 
    we see that  $\Res S_{k}'=\Res S_{k}$. Hence  $L_{k}(M,\cE)$ does not depend on the choices of $\theta$ and $J$.

   Next, let $(M',\theta')$ be a contact manifold together with a calibrated almost complex structure on $H'=\theta'$ such that there exists 
    a Heisenberg diffeomorphism $\phi$ from $(M',H')$ onto $(M,H)$. Define $\cE'=\phi^{*}\cE$ and let $S_{k}'\in \Psi_{H'}^{0}(M',\cE')$ be a 
    generalized \psivdo\ projection at level $k$ with respect $\theta'$ and a calibrated almost complex structure $J'$ on $H'$. Since $\Res S_{k}'$ depends 
    neither on $\theta$, nor on $J$, we may assume that we have $\theta'=\phi^{*}\theta$ and $J'=\phi^{*}J$. 
    
    By the results of~\cite[Sect.~3.2]{Po:MAMS1} the operator $\phi^{*}S_{k}$ is a projection in $\Psi_{H'}^{0}(M',\cE')$ with principal symbol 
    $\phi^{*}s_{k}(x,\xi)\otimes \op{id}_{\cE'}$, where $\phi^{*}s_{k}(x,\xi)=s_{k}(\phi(x),(\phi_{H}'(x)^{-1})^{t}\xi)$ and $\phi_{H}'$ is the 
    the vector bundle isomorphism  from $\fg M'=(TM'/H')\oplus H'$ onto  $\fg M=(TM/H)\oplus H$ induced by $\phi'$. 
   
   Let $s_{k}'$ be the Szeg\"o symbol on $M'$ at level $k$ with respect to $\theta'$ and $J'$. We claim that $s_{k}'=\phi^{*}s_{k}$. To see this 
 let $X_{1},\ldots,X_{2n}$ be an admissible orthonormal frame near a point $a\in M$ and for $j=0,..,2n$ let $X_{j}'=\phi^{*}X_{j}$.  Then 
        $X_{0}'$ is the  Reeb vector field of $\theta'$ and $X_{1}',\ldots,X_{2n}'$ is an admissible orthonormal frame of $H'$ near $a'=\phi^{-1}(a)$ 
        with respect to $g_{\theta',J'}=\phi^{*}g_{\theta,J}$. 
 
        Moreover, by the results of~\cite{Po:Pacific1}  for $j=0,1,\ldots,2n$ we have $(X_{j}')^{a'}=\phi_{H}'(a)^{*}X_{j}^{a}$.  Therefore,  
     $\phi_{X',a'}=\phi_{X,a}\circ \phi_{H}'(a)$ and so 
    $\phi_{X',a'}^{*}s_{k^{0}}=\phi_{H}'(a)^{*}\phi_{X,a}s_{k}^{0}=(\phi^{*}s_{k})(a,.)$. Hence $s_{k}'=\phi^{*}s_{k}$ as claimed.
   
    It follows from this that $\phi^{*}S_{k}$ is a generalized Szeg\"o projection at level $k$ on $M'$ 
    with respect to $\theta'$ and $J'$, so by the first part of the proof $\Res S_{k}'=\Res \phi_{*}S_{k}$. Since $\Res \phi_{*}S_{k}=\Res S_{k}$, 
    we see that $\Res S_{k}'=\Res S_{k}$. Hence $L_{k}(M,\cE)$ is a Heisenberg diffeomorphism invariant of $M$. 
%
%
%
%
%
   
    Let us now prove that $L_{k}(M,\cE)$ is an invariant of the $K$-theory class of $\cE$. Let $\phi$ be a vector bundle 
    isomorphism from $\cE$ onto a vector bundle $\cE'$ over $M$ and let $S_{k}'\in\pvdo^{0}(M,\cE')$ be a generalized Szeg\"o projection of level 
    $k$. Then $\phi_{*}S_{k}$ is a projection in $\pvdo^{0}(M,\cE')$ with principal symbol $s_{k}\otimes \op{id}_{\cE'}$, hence is 
    a generalized Szeg\"o projection of level $k$. Thus $\Res S_{k}'=\Res \phi_{*}S_{k}=\Res S_{k}$. 
    
    Next, for $j=1,2$ let $\cE_{j}$ be a vector bundle over $M$ and let $S_{k,\cE_{j}}\in \pvdo^{0}(M,\cE_{j})$ be a generalized Szeg\"o projection at level $k$ 
    acting on the section of $\cE_{j}$. In addition, let $S_{k,\cE_{1}\oplus\cE_{2}}\in  \pvdo^{0}(M,\cE_{1}\oplus \cE_{2})$ be a 
    generalized Szeg\"o projection at level $k$ acting on the section of $\cE_{1}\oplus \cE_{1}$. Then 
    $S_{k,\cE_{1}}\oplus S_{k,\cE_{2}}$ is a \psivdo\ projection acting on the sections of $\cE_{1}\oplus \cE_{2}$ with principal symbol $s_{k}\otimes 
    \op{id}_{\cE_{1}\oplus \cE_{2}}$, hence is a generalized Szeg\"o projection at level $k$. Thus
$\Res S_{k,\cE_{1}\oplus \cE_{2}}=\Res (S_{k,\cE_{1}}\oplus S_{k,\cE_{2}})=\Res S_{k,\cE_{1}}+\Res S_{k,\cE_{2}}$, i.e., we have 
$L_{k}(M,\cE_{1}\oplus \cE_{2})= 
L_{k}(M,\cE_{1}) + L_{k}(M,\cE_{2})$.    

    Bearing this is in mind, let $\cE'$ be a (smooth) vector bundle in the $K$-theory class of $\cE$, so that there exist an auxiliary vector bundle $\cF$ and a 
    vector bundle isomorphism $\phi$ from $\cE\oplus \cF$ onto $\cE'\oplus \cF$. Then we have $L_{k}(M, \cE\oplus \cF)=L_{k}(M, \cE'\oplus \cF)$. As 
    $L_{k}(M, \cE\oplus \cF)=L_{k}(M, \cE)+ L_{k}(M, \cF)$ and $L_{k}(M, \cE'\oplus \cF)=L_{k}(M, \cE')+ L_{k}(M, \cF)$, it follows that 
    $L_{k}(M, \cE)=L_{k}(M, \cE')$. Hence $L_{k}(M,\cE)$  is an invariant of the $K$-theory class of $\cE$.
    
    Finally, since $L_{k}(M,\cE)$ is a Heisenberg diffeomorphism invariant of $M$ arguing as in the proof of Proposition~\ref{prop:CR.deformation-spc}  
    shows its invariance under deformation of 
    the contact structure of $M$.
\end{proof}

\begin{remark}\label{rem:Szego.computation}
The almost complex structure $J$ and the Reeb vector field $X_{0}$ give rise to splittings as in~(\ref{eq:CR-decomposition}) 
and~(\ref{eq:CR-Lambda-pq-decomposition}), so that $(p,q)$-forms make sense. We then can define  generalized Szeg\"o projections on $(0,q)$-forms with 
coefficients in a vector bundle $\cE$. This can 
be done at any integer level $k=0,1,\ldots$, but for $1\leq q\leq n-1$ and $k\leq 2q-1$  the corresponding Szeg\"o symbol vanishes and we get a smoothing projection 
with vanishing noncommutative residue. Then arguing as in the proof of Theorem~\ref{thm:IGSP.main} shows that the corresponding noncommutative residues don't 
depend on the choice of the operator and yield contact invariants. 

If $(M,H)$ is CR, i.e., if $J$ defines a complex structure on $H$, then $M$ is strictly pseudoconvex, 
the contact form $\theta$ defines a pseudohermitian structure and $g_{\theta,J}$ is the associated Levi metric. Let $\cE$ be a Hermitian CR vector 
bundle equipped with a compatible CR connection $\nabla$ and let $\nabla^{0,q}_{X_{0}}$ denote the covariant derivative 
$\cL_{X_{0}}\otimes 1_{\cE}+1_{\Lambda^{0,q}}\otimes \nabla_{X_{0}}$. 
We cannot make use of the formulas~(\ref{eq:Szego-projections-dbarb-varthetab})--(\ref{eq:projections-dbarb-varthetab}) 
to prove that the projections $\Pi_{0}(\Box_{b,\cE;q}+ik\nabla^{0,q}_{X_{0}})$ are \psivdos, but the 
arguments of~\cite[\S 25]{BG:CHM} can be extended to prove this result~(\cite{Gr:PC}). Therefore, we get higher level versions of the invariants from 
Theorem~\ref{CR:Thm2}. Furthermore, in this case it should be possible to apply the Getzler's rescaling techniques alluded to in 
Subsection~\ref{subsec:CR.computation} to similarly compute the supersymmetric densities, 
\begin{equation}
    \Str_{\Lambda^{0,*}(\cE)} c_{\Pi_{0}(\Box_{b,\cE}+ik\nabla^{0,q}_{X_{0}})}(x)= 
    \sum_{q=0}^{n}(-1)^{q}\Tr_{\Lambda^{0,q}(\cE)}c_{\Pi_{0}(\Box_{b,\cE;q}+ik\nabla^{0,q}_{X_{0}})}(x).  
\end{equation}
We hope to be able to deal with this computation in a subsequent paper.
  \end{remark}

%

\section{Invariants from the contact complex}\label{sec.contact-complex}
Let $(M^{2n+1},H)$ be an orientable contact manifold. Let $\theta$ be a contact form on $M$ and let $X_{0}$ be its 
Reeb vector field of $\theta$. We also let $J$ be a calibrated almost complex structure on $H$ and  as in the previous section 
we endow $TM$ with the Riemannian metric $g_{\theta,J}=d\theta(.,J.)+\theta^{2}$.

Observe that the splitting $TM=H\oplus \R X_{0}$ allows us to identify  
$H^{*}$ with the annihilator of $X_{0}$ in $T^{*}M$. More generally, identifying $\Lambda^{k}_{\C}H^{*}$ with $\ker \iota_{X_{0}}$, where 
$\iota_{X_{0}}$ denotes the contraction operator by $X_{0}$, gives the splitting
\begin{equation}
    \Lambda^{*}_{\C}TM=(\bigoplus_{k=0}^{2n}\Lambda^{k}_{\C}H^{*}) \oplus (\bigoplus_{k=0}^{2n} \theta\wedge \Lambda^{k}_{\C}H^{*}).
     \label{eq:contact.decomposition-forms}
\end{equation}
 
For any horizontal form $\eta\in C^{\infty}(M,\Lambda^{k}_{\C}H^{*})$ we can write
$d\eta= d_{b}\eta+\theta \wedge \cL_{X_{0}}\eta$,
where $d_{b}\eta$ is the component of $d\eta$ in $\Lambda^{k}_{\C}H^{*}$. This does not provide us with a complex, for we have 
$d_{b}^{2}=-\cL_{X_{0}}\varepsilon(d\theta)=-\varepsilon (d\theta)\cL_{X_{0}}$ where $\varepsilon(d\theta)$ denotes the exterior multiplication 
by $d\theta$.

%

The contact complex of Rumin~\cite{Ru:FDVC} is an attempt to get 
a complex of horizontal differential forms by forcing the equalities $d_{b}^{2}=0$ and $(d^{*}_{b})^{2}=0$.

A natural way to modify $d_{b}$ to get the equality $d_{b}^{2}=0$ is to restrict 
$d_{b}$ to the subbundle $\Lambda^{*}_{2}:=\ker \varepsilon(d\theta) \cap \Lambda^{*}_{\C}H^{*}$, since the latter
is closed under $d_{b}$ and is annihilated by $d_{b}^{2}$. 

Similarly, we get the equality $(d_{b}^{*})^{2}=0$ by restricting $d^{*}_{b}$ to the subbundle 
$\Lambda^{*}_{1}:=\ker \iota(d\theta)\cap \Lambda^{*}_{\C}H^{*}=(\im \varepsilon(d\theta))^{\perp}\cap \Lambda^{*}_{\C}H^{*}$, where 
$\iota(d\theta)$ denotes the interior product 
with $d\theta$. This amounts to replace $d_{b}$ by $\pi_{1}\circ d_{b}$, where $\pi_{1}$ is the orthogonal projection onto $\Lambda^{*}_{1}$.

In fact, since $d\theta$ is nondegenerate on $H$ the operator $\varepsilon(d\theta):\Lambda^{k}_{\C}H^{*}\rightarrow \Lambda^{k+2}_{\C}H^{*}$  is 
injective for $k\leq n-1$ and surjective for $k\geq n+1$. This implies that $\Lambda_{2}^{k}=0$ for $k\leq n$ and $\Lambda_{1}^{k}=0$ for $k\geq n+1$. 
Therefore, we only have two halves of complexes. 

As observed by Rumin~\cite{Ru:FDVC} we get a full complex by connecting 
the two halves by means of the operator $D_{R,n}:C^{\infty}(M,\Lambda_{\C}^{n}H^{*}) \rightarrow 
C^{\infty}(M,\Lambda_{\C}^{n}H^{*})$ such that 
\begin{equation}
    D_{R,n}=\cL_{X_{0}}+d_{b,n-1}\varepsilon(d\theta)^{-1}d_{b,n},
\end{equation}
where $\varepsilon(d\theta)^{-1}$ is the inverse of $\varepsilon(d\theta):\Lambda^{n-1}_{\C}H^{*}\rightarrow \Lambda^{n+1}_{\C}H^{*}$. Notice that 
$D_{R,n}$ is second order differential 
operator.  Thus, if we let $\Lambda^{k}=\Lambda_{1}^{k}$ for $k=0,\ldots,n-1$ and 
$\Lambda^{k}=\Lambda_{1}^{k}$ for $k=n+1,\ldots,2n$ then we get the contact complex, 
\begin{multline}
    C^{\infty}(M)\stackrel{d_{R;0}}{\rightarrow}C^{\infty}(M,\Lambda^{1})\stackrel{d_{R;1}}{\rightarrow}
    \ldots 
    C^{\infty}(M,\Lambda^{n-1})\stackrel{d_{R;n-1}}{\rightarrow}C^{\infty}(M,\Lambda^{n}_{1})\stackrel{B_{R}}{\rightarrow}\\  
    C^{\infty}(M,\Lambda^{n}_{2})  \stackrel{d_{R;n}}{\rightarrow}C^{\infty}(M,\Lambda^{n+1}_{2})
    \ldots \stackrel{d_{R;2n-1}}{\longrightarrow} C^{\infty}(M,\Lambda^{2n}),
     \label{eq:contact-complex}
\end{multline}
where $d_{R;k}=\pi_{1}\circ d_{b;k}$ for $k=0,\ldots,n-1$ (so that $d_{R;k+1}^{*}=d_{b;k+1}^{*}$) and $d_{R;k}=d_{b;k}$ for $k=n,\ldots,2n-1$. 

The contact Laplacian is defined as follows. In degree $k\neq n$ this is the differential operator 
$\Delta_{R,k}:C^{\infty}(M,\Lambda^{k})\rightarrow C^{\infty}(M,\Lambda^{k})$ such that
\begin{equation}
    \Delta_{R,k}=\left\{
    \begin{array}{ll}
        (n-k)d_{R,k-1}d^{*}_{R,k}+(n-k+1) d^{*}_{R,k+1}d_{R,k}& \text{$k=0,\ldots,n-1$},\\
         (k-n-1)d_{R,k-1}d^{*}_{R,k}+(k-n) d^{*}_{R,k+1}d_{R,k}& \text{$k=n+1,\ldots,2n$}.
         \label{eq:contact-Laplacian1}
    \end{array}\right.
\end{equation}
For $k=n$ we have the differential operators $\Delta_{R,nj}:C^{\infty}(M,\Lambda_{j}^{n})\rightarrow C^{\infty}(M,\Lambda^{n}_{j})$, $j=1,2$, 
given by the formulas, 
\begin{equation}
    \Delta_{R,n1}= (d_{R,n-1}d^{*}_{R,n})^{2}+D_{R,n}^{*}D_{R,n}, \quad   \Delta_{R,n2}=D_{R,n}D_{R,n}^{*}+  (d^{*}_{R,n+1}d_{R,n}).
         \label{eq:contact-Laplacian2}
\end{equation}

Observe that $\Delta_{R,k}$, $k\neq n$, is a differential operator order $2$, whereas $\Delta_{Rn1}$ and $\Delta_{Rn2}$ are differential operators of 
order $4$. Moreover, Rumin~\cite{Ru:FDVC} proved that in every degree the contact Laplacian is maximal hypoelliptic. 
In fact, in every degree the contact Laplacian has an invertible principal symbol, hence admits a parametrix in the Heisenberg calculus 
(see~\cite{JK:OKTGSU},~\cite[Sect.~3.5]{Po:MAMS1}).
 
Let $\Pi_{0}(d_{R,k})$ and $\Pi_{0}(D_{R,n})$ be the orthogonal projections onto $\ker d_{R,k}$ and $\ker D_{R,n}$, and let 
$\Delta_{R,k}^{-1}$ and $\Delta_{R,nj}^{-1}$ be the partial inverses of $\Delta_{R,k}$ and $\Delta_{R,nj}$.
Then as in~(\ref{eq:projections-dbarb-varthetab}) we have 
%
 \begin{gather}
     \Pi_{0}(d_{R,k})=\left\{ 
     \begin{array}{ll}
         1-(n-k-1)^{-1}d^{*}_{R,k+1}\Delta_{R,k+1}^{-1}d_{R,k} & k=0,\ldots,n-2,\\
          1- d^{*}_{R,n}d_{R,n-1}d^{*}_{R,n}\Delta_{R,n1}^{-1}d_{R,n-1} & k=n-1,\\
          1- (k-n)^{-1}d^{*}_{R,k+1}\Delta_{R,k+1}^{-1}d_{R,k} & \text{$k=n,\ldots,2n-1$},
     \end{array} \right. \\
     \Pi_{0}(D_{R,n})=1-D_{R,*}\Delta_{R,n2}^{-1}D_{R,n}.
  \end{gather}

As in each degree the principal symbol of the contact Laplacian is invertible, 
the operators $\Delta_{R,k}^{-1}$, $k\neq n$, and $\Delta_{R,nj}^{-1}$, $j=1,2$  are \psivdos\ of order $-2$ and order $-4$ respectively. 
Therefore, the above formulas for $\Pi_{0}(d_{R,k})$ and $\Pi_{0}(D_{R,n})$ show that these projections are zero'th order \psivdos.

\begin{theorem}\label{Contact:Thm}
    The noncommutative residues $\Res \Pi_{0}(d_{R,k})$, $k=1,\ldots,2n-1$ and $\Res \Pi_{0}(D_{R,n})$ are Heisenberg diffeomorphism invariants of 
    $M$ and are invariant under deformation of the contact structure. In particular, their values depend neither on the contact form $\theta$, nor on the almost 
    complex structure $J$. 
\end{theorem}
\begin{proof}
    Let us first show that the above noncommutative residues don't depend 
    on $\theta$ or on $J$. Let $\theta'$ be a contact form on $M$ and let $J'$ be an almost complex structure on $H$ calibrated with respect to $\theta'$. 
    We shall assign the superscript $'$ to objects associated to the pair $(\theta',J')$. 
    
    It follows from~\cite{Ru:FDVC} and \cite[Chap.~2]{Po:MAMS1} that there are vector bundle isomorphisms 
    $\phi_{k}:\Lambda^{k}_{1}\rightarrow\Lambda^{k'}_{1}$, $k=1,\ldots, n$, and $\psi_{l}:\Lambda^{l}_{2}\rightarrow \Lambda^{l'}_{2}$, 
    $l=n,\ldots,2n-1$
    such that $\phi_{k+1}d_{R,k}=d_{R,k}'\phi_{k}$ and $\psi_{l+1}d_{R,l}=d_{R,l}'\psi_{l}$ and also $\psi_{n}D_{R,n}=D_{R,n}'\phi_{n}$. 
    In particular, we have $\ker d_{R,k}'=\phi_{k}(\ker d_{R,k})$. Since $\ker d_{R,k}$ is the range 
   of $\Pi_{0}(d_{R,k}')$ and $\phi_{k}(\ker d_{R,k})$ that of $\phi_{k*}\Pi_{0}(d_{R,k})$, we see that $\Pi_{0}(d_{R,k}')$ and 
   $\phi_{k*}\Pi_{0}(d_{R,k})$ have same range, so by Proposition~\ref{prop:invariance-range-kernel} 
     $  \Res \Pi_{0}(d_{R,k}')=\Res \phi_{k*}\Pi_{0}(d_{R,k})=\Res \Pi_{0}(d_{R,k})$.
     
          Similarly, $\Res \Pi_{0}(d_{R,l}')=\Res \Pi_{0}(d_{R,l})$ and 
          $\Res \Pi_{0}(D_{R,n}')=\Res \Pi_{0}(D_{R,n})$. Hence
     $\Res \Pi_{0}(d_{R,k})$, $k=1,\ldots,2n-1$, and $\Res \Pi_{0}(D_{R,n})$ 
     don't depend on $\theta$ or $J$. 

Next,  let $(M',H')$ be an orientable contact manifold and let $\phi:M'\rightarrow M$ be  a Heisenberg diffeomorphism.   
In addition, let $\theta'$ be a contact form on $M'$ and let $J'$ be a calibrated almost complex structure on $H'$. 
As before we shall assign the superscript $'$ to objects related to $M'$. 

Since 
$\Res \Pi_{0}(d_{R,k})$, $k=1,\ldots,2n-1$, and $\Res \Pi_{0}(D_{R,n})$  
don't depend on the contact form or on the almost complex structure, we may assume that $\theta'=\phi^{*}\theta$ and 
$J'=\phi^{*}J$. Then $\phi^{*}g_{\theta,J}=g_{\theta',J'}$ and $\phi^{*}d_{R,k}=d_{R,k}'$. Thus, 
\begin{equation}
    \phi^{*}\Pi_{0}(D_{R,n})=\Pi_{0}(D_{R,n}'), \qquad \phi^{*}\Pi_{0}(d_{R,k})=\Pi_{0}(d_{R,k}'). 
\end{equation}
Since $\Res \phi^{*}\Pi_{0}(d_{R,k})=\Res \Pi_{0}(d_{R,k})$ it follows that $\Res \Pi_{0}(d_{R,k}')=\Res \Pi_{0}(d_{R,k})$. Similarly, we have 
$\Res \phi^{*}\Pi_{0}(D_{R,n}')=\Res \Pi_{0}(D_{R,n})$.
Hence  $\Res \Pi_{0}(D_{R,n})$ and $\Res \Pi_{0}(d_{R,k})$ are Heisenberg diffeomorphism invariants.  

Finally, since the noncommutative residues  $\Res \Pi_{0}(D_{R,n})$ and $\Res \Pi_{0}(d_{R,k})$ are invariant by Heisenberg diffeomorphism we may argue as in the 
proof  of Proposition~\ref{prop:CR.deformation-spc} to show that they are invariant under deformation of the contact structure. The proof is thus achieved.
\end{proof}

 \begin{remark}
The residues $\Res \Pi_{0}(d_{R,k}^{*})$, $k=1,\ldots,2n-1$, and $\Res \Pi_{0}(D_{R,n}^{*})$ also yield invariants, but these are the same up to a 
sign factor as those from Theorem~\ref{Contact:Thm}. 
Indeed, in~(\ref{eq:Szego-projections-dbarb-varthetab}) for $k\neq n$ we have  $\Pi_{0}(\Delta_{R,k})=\Pi_{0}(d_{R,k})+\Pi_{0}(d^{*}_{R,k})-1$, 
so as $\Pi_{0}(\Delta_{R,k})$ is a  
 smoothing operator we get $\Res \Pi_{0}(d^{*}_{R,k})=-\Res \Pi_{0}(d_{R,k})$. Similarly, we have $ \Res \Pi_{0}(d^{*}_{R,n})=- \Res  \Pi_{0}(D_{R,n})$  and 
 $ \Res \Pi_{0}(D_{R,n}^{*})= -\Res \Pi_{0}(d_{R,n})$. 
\end{remark} 

\begin{remark}
    As with the invariants of the previous sections, we can make use of a Getzler's rescaling to compute the local densities associated to a 
    supersymmtric version of the invariants of Theorem~\ref{Contact:Thm}. However, it is not clear how efficient this would be to yield explicit formulas. 
    More precisely, in the non-supersymetric setting 
there are no known explicit formulas for the fundamental solutions of the contact Laplacian and \emph{a fortiori} for the inverse of its principal 
symbols. Therefore, it is all the more difficult to obtain explicit formulas in the supersymmetric setting. 

Maybe the solution would be combine a Getzler's rescaling with diabatic limit techniques, since the contact complex appears naturally in the analysis of the 
asymptotic behavior of the de Rham complex under a diabatic limit (see~\cite{Ru:SRLDFSCM}, \cite{BHR:DLEICRMD3}). Let us also mention that another 
possible approach would be rely on the global $K$-theoretic techniques as alluded to in Appendix.
\end{remark}

\appendix
\setcounter{section}{1}
\section*{Appendix}
\label{sec:K-theory}
The construction of the contact invariants from generalized Szeg\"o projection in Section~\ref{sec.Toeplitz} is partly based on 
Lemma~\ref{lem:NCRPP.homotopy-invariance} 
stating that the 
noncommutative residue of a \psivdo\ projection is a homotopy invariant of its principal symbol. 
In this Appendix we would like to explain that for a general Heisenberg manifold, not necessarily contact or CR, this leads us to 
 a $K$-theoretic interpretation of the noncommutative residue of a \psivdo\ projection. 

 First, as $S_{0}(\fg^{*}M)$  endowed with the 
 product homogeneous symbols~(\ref{eq:NCRP.product-symbols}) is a Fr\'echet algebra, its $K_{0}$-group can be defined as follows.
 
 Let $M_{\infty}(S_{0}(\fg^{*}M))$ be the algebra $\varinjlim M_{q}(S_{0}(\fg^{*}M))$, where the inductive limit is defined using the 
 embedding $ a \rightarrow \op{diag}(a,0)$ of 
 $M_{q}(S_{0}(\fg^{*}M))$ into $M_{q+1}(S_{0}(\fg^{*}M))$
 We say that idempotents $e_{1}$ and $e_{2}$ in $M_{\infty}(S_{0}(\fg^{*}M))$  are homotopic 
  when there is a $C^{1}$-path of idempotents joining $e_{1}$ to $e_{2}$ is some space $M_{q}(S_{0}(\fg^{*}M))$ containing both 
  $e_{1}$ and $e_{2}$. 
 
 The addition of idempotents is given by the direct sum \(e_{1}\oplus e_{2}=\op{diag}(e_{1},e_{2})\).  
This operation is compatible with the homotopy of idempotents, so turns the set of homotopy classes of  
idempotents in $M_{\infty}(S_{0}(\fg^{*}M))$ into a monoid. We then define $K_{0}(S_{0}(\fg^{*}M))$ as the associated 
Abelian group of this monoid, i.e., the group of formal 
differences of homotopy classes of idempotents in $M_{\infty}(S_{0}(\fg^{*}M))$. 
 
In this sequel we will need to describe $K_{0}(S_{0}(\fg^{*}M))$ in terms of homotopy classes of symbols as follows. 
 
\begin{definition}
1) An idempotent pair is a pair $(\pi,\cE)$ consisting of a (smooth) vector bundle $\cE$ over $M$ and an idempotent $\pi \in S_{0}(\fg^{*}M,\cE)$.\smallskip

2) Two idempotent pairs $(\pi_{1},\cE_{1})$ and $(\pi_{2},\cE_{2})$ are equivalent when there exist smooth vector bundles $\cF_{1}$ and $\cF_{2}$ 
over $M$, a vector bundle isomorphism $\phi$ from $\cE_{1}\oplus \cF_{1}$ to $\cE_{2}\oplus \cF_{2}$ and a homotopy of idempotents in 
$S_{0}(\fg^{*}M,\cE_{2}\oplus \cF_{2})$ from $\phi_{*}(\pi_{1}\oplus 0_{\pr^{*} \cF_{1}})$ to $\pi_{2}\oplus 0_{\pr^{*} \cF_{2}}$.
\end{definition}
 

The  set of equivalence classes of idempotent pairs becomes a monoid when we endow it with the addition given by the direct sum, 
 \begin{equation}
     (\pi_{1},\cE_{1})\oplus (\pi_{2},\cE_{2})=(\pi_{1}\oplus\pi_{2},\cE_{1}\oplus\cE_{2}).
 \end{equation}
We let $\cI_{0}(\fg^{*}M)$ denote the Abelian group of formal differences of equivalence classes of idempotent pairs. 

There is a natural map $\Theta$ from idempotents of $M_{q}(S_{0}(\fg^{*}M))$ to idempotent pairs obtained by assigning to any idempotent 
$e \in M_{q}(S_{0}(\fg_{*}M))$ the idempotent pair $\Theta(e)=(e,M\times \C^{q})$. In fact, we have:

\begin{lemma}
   $\Theta$ gives rise to an isomorphism from $K_{0}(S_{0}(\fg^{*}M))$ onto $\cI_{0}(\fg^{*}M)$.
\end{lemma}
\begin{proof}
    First, let $(\pi,\cE)$ be an idempotent pair. There  exists a (smooth) vector bundle $\cF$ such 
    that $\cE\oplus \cF$ is globally trivializable, i.e., there exists a vector bundle isomorphism $\phi:\cE\oplus \cF\simeq 
 M\times \C^{q}$ (see, e.g.,~\cite[Cor.~1.4.14]{At:KT}). Then $(\pi,\cE)$ is equivalent to $\Theta(e)=(e,M\times \C^{q})$, where $e$ is the idempotent 
 $\phi_{*}(\pi\oplus 0_{\pr \cF})\in M_{q}(S_{0}(\fg^{*}M))$. This shows that up to the equivalence of idempotent pairs the map $\Theta$ is surjective. 
 
 Next, for $j=1,2$ let $e_{j}\in M_{q_{j}}(S_{0}(\fg^{*}M))$ be idempotent and suppose that $\Theta(e_{1})$ and $\Theta(e_{2})$ are equivalent 
 idempotent pairs. Thus there exist vector bundles $\cE_{1}$ and $\cE_{2}$ 
 and a vector bundle isomorphism $\phi:(M\times \C^{q_{1}})\oplus 
 \cE_{1}\rightarrow (M\times \C^{q_{2}})\oplus \cE_{2}$ such that $\phi_{*}(e_{1}\oplus 0_{\pr^{*}\cE_{1}})$ is homotopic to $e_{2}\oplus 
 0_{\pr^{*}\cE_{2}}$ in $S_{0}(\fg^{*}M,(M\times \C^{q_{2}})\oplus \cE_{2})$. 
 
 Let $\cF$ be a vector bundle so that there exists a vector bundle 
 isomorphism $\psi_{1}$ from $(M\times \C^{q_{1}})\oplus \cE_{1}\oplus \cF$ onto $M\times \C^{q}$. Composing $\psi_{1}$ 
 with $\phi\oplus 1_{\cF}$ we get 
 a vector bundle isomorphism $\psi_{2}: (M\times \C^{q_{2}})\oplus \cE_{2}\oplus \cF \rightarrow M\times \C^{q}$ such that there is a homotopy of 
 idempotents in $S_{0}(\fg^{*}, M\times \C^{q})=M_{q}(S_{0}(\fg^{*}))$ joining $\psi_{1*}(e_{1}\oplus 
 0_{\pr^{*}(\cE_{1}\oplus \cF)})$ to  $\psi_{2*}(e_{2}\oplus  0_{\pr^{*}(\cE_{2}\oplus \cF)})$.
 
Using the identification $\C^{q}=\C^{q_{1}}\oplus \C^{q-q_{1}}$ let  $\gamma=1_{\C^{q}}\oplus 0$, where $0$ is the zero vector bundle morphism from 
$\cE_{1}\oplus \cF$ to $M\times \C^{q-q_{1}}$, and similarly define $\gamma^{t}= 1_{\C^{q}}\oplus 0$ using the the zero vector bundle morphism from 
$M\times \C^{q-q_{1}}$ to $\cE_{1}\oplus \cF$. Let $\alpha=\psi_{1}\circ \gamma$ and $\beta=\gamma^{t}\circ \psi^{-1}_{1}$. Then $\alpha$ and $\beta$ 
are sections of $\End (M\times \C^{q})$, i.e., are elements of $M_{q}(C^{\infty}(M))$, and we have 
$\alpha(e_{1}\oplus 0_{\C^{q-q_{1}}})\beta=\psi_{1*}(e_{1}\oplus  0_{\pr^{*}(\cE_{1}\oplus \cF)})$ and $(e_{1}\oplus  0_{\C^{q-q_{1}}})\beta 
\alpha=e_{1}\oplus  0_{\C^{q-q_{1}}}$. This shows  that $e_{1}\oplus 0_{\C^{q-q_{1}}}$ and $\psi_{1*}(e_{1}\oplus  0_{\pr^{*}(\cE_{1}\oplus \cF)})$ are 
algebraically equivalent idempotents of $M_{q}(S_{0}(\fg^{*}M))$ in the sense of~\cite[Def.~4.2.1]{Bl:KTOA}. As this implies that there 
are homotopic idempotents in $M_{\infty}(S_{0}(\fg^{*}M))$ (see~\cite[Props.~4.3.1, 4.4.1]{Bl:KTOA}), it follows that
$e_{1}$ and $\psi_{1*}(e_{1}\oplus  0_{\pr^{*}(\cE_{1}\oplus \cF)})$ are homotopic to each other.

Similarly, $e_{2}$ and $\psi_{2*}(e_{2}\oplus  0_{\pr^{*}(\cE_{2}\oplus \cF)})$ are homotopic idempotents. Since the 
latter is homotopic to $\psi_{1*}(e_{1}\oplus  0_{\pr^{*}(\cE_{1}\oplus \cF)})$, it follows that $e_{1}$ is homotopic to $e_{2}$. 
This proves that $\Theta$ is injective up to the equivalence of idempotent pairs. 


All this shows that $\Theta$ factorizes through a bijection from homotopy classes of idempotents in $M_{\infty}(S_{0}(\fg^{*}M))$ to equivalence classes of 
idempotent pairs. Furthermore, this map is additive. Indeed, if $e_{j}\in M_{q_{j}}(S_{0}(\fg^{*}M))$ then, up to the identification 
$\C^{q_{1}}\oplus \C^{q_{2}}=\C^{q_{1}+q_{2}}$, we have $\Theta(e_{1}\oplus e_{2})=\Theta(e_{1})\oplus \Theta(e_{2})$. 
Thus $\Theta$ gives rise to isomorphism from $K_{0}(S_{0}(\fg^{*}M))$ onto $\cI_{0}(\fg^{*}M)$.   
\end{proof}

We will also need the following lemma. 
 \begin{lemma}\label{lem:K-Theory.symbol-projection}
     Let $\pi_{0}\in S_{0}(\fg^{*}M,\cE)$ be idempotent. Then we can always find a \psivdo\ projection $\Pi\in\pvdo^{0}(M,\cE)$ with principal symbol 
     $\pi_{0}$. 
 \end{lemma}
\begin{proof}
    Let $f_{0}=\frac{1}{2}(1+\pi_{0})$, so that $f_{0}*f_{0}=1$. Since the principal symbol map $\sigma_{0}: \pvdo^{0}(M,\cE)\rightarrow 
    S_{0}(\fg^{*}M)$ is surjective (see~\cite[Prop.~3.2.6]{Po:MAMS1}), there exists $P\in \pvdo^{0}(M,\cE)$ with principal symbol $f_{0}$. 
    Then $P^{2}$ has principal symbol 
    $f_{0}*f_{0}=1$, so $P^{2}=1-R_{1}$ with $R_{1}\in \pvdo^{-1}(M,\cE)$. 
    
    Let $\sum_{k\geq 0}a_{k}z^{k}$ be the Taylor series at $z=0$ of $(1-z)^{-\frac{1}{2}}$. Since $R_{1}$ has order~$\leq -1$ the symbolic calculus 
    for \psivdos\ in~\cite{BG:CHM} allows us to construct $Q\in \pvdo^{0}(M,\cE)$ such that 
    $Q=\sum_{k=0}^{N-1}a_{k}R^{k}_{1}\bmod \pvdo^{-N}(M,\cE)$ for any integer $N$.
     Then we obtain $(1-R_{1})Q^{2}=1 \bmod \psinf(M,\cE)$ and, letting $F=PQ$ and using the fact that $R_{1}=1-P^{2}$ commutes with $P$, we see that 
    $F^{2}=P^{2}Q=(1-R_{1})Q=1-R$ for some smoothing operator $R$.
    
    Next, as for  $\lambda\in \C$ we have $(F-\lambda)(F+\lambda)=F^{2}-\lambda^{2}=1-\lambda^{2}-R$, we see that $\lambda \in \op{Sp}F$ if, and only 
    if, $\lambda^{2}-1\in \op{Sp}R$. Since $R$ is smoothing this is a compact operator and so $\Sp R\setminus \{0\}$ is bounded and 
    discrete. Incidentally $\Sp F\setminus\{\pm 1\}$ is a discrete set. Moreover, for $\lambda \not \in (\Sp F \cup \{\pm 1\})$ we have 
    \begin{equation}
        (F-\lambda)^{-1}=(F+\lambda)(1-\lambda^{2}-R)^{2}=(\lambda^{2}-1)^{-1}(F+\lambda)-(F+\lambda)S(1-\lambda^{2}),
         \label{eq:K-theory.resolvent-F}
    \end{equation}
    where for $\mu \not \in \op{Sp} R\cup\{0\}$ we have let $S(\mu)=(\mu-R)^{-1}-\mu^{-1}$. 
    
   At first glance $(S(\mu))_{\mu \not \in \op{Sp} R\cup\{0\}}$ is a holomorphic family of bounded operators, but the equalities 
   $S(\mu)=\mu^{-1}R(\mu-R)^{-1}=\mu^{-1}(\mu-R)^{-1}R$ imply that it actually is  a holomorphic family of smoothing operators. 
    
    Now, since $\Sp F\setminus\{\pm 1\}$ is discrete we can find positive numbers $0<r_{1}<r_{2}<2$ such that the region $r_{1}<|\lambda-1|<r_{2}$ is 
    contained in the complement of $\op{Sp}F$, so that $(F-\lambda)^{-1}$ is a holomorphic family with values in $\pvdo^{0}(M,\cE)$ on that region. Therefore, for 
    $r\in (r_{1},r_{2})$ we define a projection in $L^{2}(M,\cE)$ by letting
 \begin{equation}
        \Pi=\frac{1}{2i\pi}\int_{|\lambda-1|=r}(F-\lambda)^{-1}d\lambda.
   \end{equation}
   In fact, as $(S(\mu))_{\mu \not \in \op{Sp} R\cup\{0\}}$ is a holomorphic family of smoothing 
    operators, it follows from~(\ref{eq:K-theory.resolvent-F}) that, up to a smoothing operator, $\Pi$ agrees with 
\begin{equation}
        \frac{1}{2i\pi}\int_{|\lambda-1|=r}(\lambda^{2}-1)^{-1}(F+\lambda)d\lambda=\frac{1}{2}(F+1).
\end{equation}
     Hence $\Pi$ is a zero'th order \psivdo\ projection with 
    principal symbol $\frac{1}{2}(f_{0}+1)=\pi_{0}$. The lemma is thus proved.
    \end{proof}

We can now give the topological interpretation of the noncommutative residue of a \psivdo\ projection.

\begin{proposition}\label{prop:NCR.rho-R}
  There exists a unique additive map $\rho_R :K_{0}(S_{0}(M))\rightarrow \R$ such that, for any vector bundle $\cE$ over $M$ and any 
     projection $\Pi\in \pvdo^{0}(M,\cE)$, we have 
     \begin{equation}
         \rho_R \circ \Theta^{-1}[\pi_{0},\cE]=\Res \Pi,
         \label{eq:K-Theory.residue-morphism}
     \end{equation}
      where $\pi_{0}$ denotes the principal symbol of $\Pi$.
\end{proposition}  
\begin{proof}
    Let $(\pi,\cE)$ be an idempotent pair. By Lemma~\ref{lem:K-Theory.symbol-projection} 
 there exists a projection $\Pi_{(\pi,\cE)}$ in $\pvdo^{0}(M,\cE)$ whose principal 
 symbol is $\pi$. The choice of $\Pi_{(\pi,\cE)}$ is not unique, but Lemma~\ref{lem:NCRPP.homotopy-invariance} 
 insures us that the value of $\Res \Pi_{(\pi,\cE)}$ is independent of 
 this choice. Furthermore, we know by Lemma~\ref{NCRP.real-valued}
  that $\Res  \Pi_{(\pi,\cE)}$ is a real number. 
 Therefore, we uniquely define a map from idempotent pairs to $\R$ by letting $\rho_R'(\pi,\cE)=\Res \Pi_{(\pi,\cE)}$.
 
 For $j=1,2$ let $(\pi_{j},\cE_{j})$ be an idempotent pair and let $\Pi_{j}\in \pvdo^{0}(M,\cE_{j})$ be 
a \psivdo\ projection with principal symbol $\pi_{j}$. Then $\Pi_{1}\oplus\Pi_{2}$ is a \psivdo\ projection 
 with principal symbol $\pi_{1}\oplus \pi_{2}$, so 
  $ \rho_R'(\pi_{1}\oplus\pi_{2},\cE_{1}\oplus \cE_{2})=\Res (\Pi_{1}\oplus \Pi_{2})$. As we have 
  $\Res (\Pi_{1}\oplus \Pi_{2})=\Res 
       \Pi_{1}+\Res \Pi_{2}$ we get $\rho_R'(\pi_{1}\oplus\pi_{2},\cE_{1}\oplus \cE_{2})=\rho_R'(\pi_{1},\cE_{1})+\rho_R'(\pi_{2},\cE_{2})$.
 Hence $\rho_R'$ is an additive map. 

Next, assume that $(\pi_{1},\cE_{1})$ and $(\pi_{2},\cE_{2})$ are equivalent idempotent pairs. Thus, there exist smooth 
      vector bundles $\cF_{1}$ and $\cF_{2}$ over $M$ and a vector bundle isomorphism $\phi:\cE_{1}\oplus \cF_{1}\rightarrow
      \cE_{2}\oplus \cF_{2}$ such that $\phi_{*}(\pi_{1}\oplus 0_{\pr^{*} \cF_{1}})$ and $\pi_{2}\oplus 0_{\pr^{*} \cF_{2}}$
      are homotopic idempotents in $S_{0}(\fg^{*}M,\cE_{2}\oplus \cF_{2})$. Then 
      by Lemma~\ref{lem:NCRPP.homotopy-invariance} we have 
      $\rho_R'(\pi_{1}\oplus 0_{\pr^{*} \cF_{1}},\cE_{1}\oplus \cF_{1})=\rho_R'(\pi_{2}\oplus 0_{\pr^{*} \cF_{2}},\cE_{2}\oplus \cF_{2})$ and, as 
      $\rho_R'(0_{\pr^{*} \cF_{j}},\cF_{j})=0$, it follows from the additivity of $\rho_R'$ that $\rho_R'(\pi_{j}\oplus 0_{\pr^{*} \cF_{j}},\cE_{j}\oplus 
      \cF_{j})=\rho_R'(\pi_{j},\cE_{j})$. Hence we have $\rho_R'(\pi_{1},\cE_{1})=\rho_R'(\pi_{2},\cE_{2})$. 
      
      This shows that the value of $\rho_R'(\pi_{1},\cE_{1})$ depends 
      only on the equivalence class of $(\pi_{1},\cE_{1})$. Since $\rho_R'$ is additive it follows that it gives rise to an additive map
from $\cI_{0}(\fg^{*}M)$ to $\R$.  
      Letting $\rho_R=\rho_R'\circ\Theta$ then defines the desired additive map from $K_{0}(S_{0}(\fg^{*}M))$ to $\R$ 
      satisfying~(\ref{eq:K-Theory.residue-morphism}).
%
\end{proof}

The above $K$-theoretic interpretation of the noncommutative residue of a \psivdo\ projection is reminiscent of the $K$-theoretic interpretations of the 
residue at the origin of the eta function of a selfadjoint elliptic \psido\ by Atiyah-Patodi-Singer~\cite{APS:SARG3} and of the Fredholm 
index of an elliptic \psido\ by Atiyah-Singer~\cite{AS:IEO1}. Nevertheless, it differs from them 
on the fact that we have to use the $K$-theory of algebras rather than 
that of spaces as in~\cite{AS:IEO1} and~\cite{APS:SARG3}. 
Indeed, as the algebra of (scalar) zero'th order Heisenberg symbols is not 
commutative, it cannot be identified with the algebra of smooth functions on the cotangent unit sphere $S^{*}M$ 
and we cannot make use of the Serre-Swan isomorphism to identify its $K$-theory with that of $S^{*}M$. Thus in order to give a 
$K$-theoretic interpretation of the noncommutative residue of a \psido\ projection we really have to rely on the 
$K$-theory of algebras. 

The (full) index theorem of Atiyah-Singer~\cite{AS:IEO1} identifies in purely topological terms the $K$-theoretic analytical index map 
defined the Fredholm indices 
of elliptic \psidos. In turn, via a cohomological  interpretation this provides us with a general topological formula to compute the index of an elliptic \psido\ in 
terms of the Chern character of its principal symbol. 

Similarly, it would be interesting to have a topological formula for computing the noncommutative residue of a 
\psivdo\ projection in terms of its principal symbol. As above-mentioned the algebra of zero'th order Heisenberg symbols is noncommutative, 
so we presumably have to rely on tools from Connes' noncommutative geometry to carry out this project. In particular, Connes~\cite{Co:NCG} 
produced a fairly simple and general proof of the Atiyah-Singer index theorem by making use of the 
tangent groupoid of a manifold. The latter construction has been extended to Heisenberg manifolds in~\cite{Po:Pacific1} (see also~\cite{Va:PhD}). 
Therefore, the tangent groupoid of a Heisenberg manifold may well be a key tool to give a topological 
interpretation of the residue map $\rho_R$. 

On the other hand, by a celebrated result of Atiyah-Patodi-Singer~\cite{APS:SARG3} and Gilkey~\cite{Gi:RGEFO} the eta function of a general 
selfadjoint elliptic \psido\ on a compact manifold is regular at the origin, so that the eta invariant of the operator is always well defined. This 
result was extended by 
Wodzicki~\cite{Wo:LISA} who established the vanishing of the noncommutative residue of a \psido\ projection. 
The original proof of Wodzicki is quite involved, but it was much simplified by Br\"uening-Lesch~\cite[Lem.~2.7]{BL:OEICNLBVP} 
who showed that the result of Wodzicki is in fact equivalent to that of Atiyah-Patodi-Singer and Gilkey. 

Similarly, in the framework of the Heisenberg calculus the vanishing of the noncommutative residue of a \psivdo\ projection is equivalent to the 
regularity at the origin of the eta function of a selfadjoint hypoelliptic \psivdo. Therefore, proving the vanishing of the map $\rho_R$ would enable us to define the eta 
invariant of a selfadjoint  hypoelliptic \psivdo\ as the regular value at the origin of its eta function. Such a result would be interesting 
for dealing with hypoelliptic boundary values index problems on bounded strictly pseudoconvex complex domains, symplectic manifolds or 
even asymptotically complex hyperbolic spaces in the sense of~\cite{EMM:RLSPD}. This would also allow us to give a positive answer to a question left 
open in~\cite[Remark~9.3]{BHR:DLEICRMD3}.

To summarize two interesting phenomena may occur: 

- The map $\rho_R$ may be non-trivial and understood in topological terms, which would allows us to compute the noncommutative residue of \psivdo\ in 
terms of its principal symbol only;

- The noncommutative residue of a \psivdo\ is always zero vanish, which would allow us to define the eta invariant of any hypoelliptic selfadjoint \psidos. 

\noindent Therefore, it is all the more important to further understand the noncommutative residue of a \psivdo\ projection. We hope to go back to this in a 
future research.


\begin{acknowledgements} 
I am indebted to Louis Boutet de Monvel, Charlie Epstein, Charlie Fefferman, Peter Greiner, Kengo Hirachi and C.~Robin Graham for 
stimulating discussions about the subject matter of the paper.  I would also like to thank for their hospitality the University of California at Berkeley and the University 
of Tokyo where parts of the paper were written. 
\end{acknowledgements}

\end{document}